\documentclass[11pt,a4paper]{article}
\usepackage[utf8]{inputenc}
\usepackage{amsmath,bm}
\usepackage{amsmath}
\usepackage{amsthm}
\usepackage{amsfonts}
\usepackage[square, comma, sort&compress,numbers]{natbib}
\usepackage{amssymb}
\usepackage{setspace}
\usepackage{graphicx}
\usepackage{verbatim} 
\usepackage{color}
\usepackage[marginal]{footmisc}

\usepackage{indentfirst}
\usepackage[left=2cm,right=2cm,top=2cm,bottom=2cm]{geometry}
\author{}
\usepackage{float}
\newtheorem{lem}{Lemma}
\newtheorem{thm}{Theorem}

\begin{document}
\bibliographystyle{plainnat}

\title{\textbf {Homoclinic Bifurcations that Give Rise to Heterodimensional Cycles near A Saddle-Focus Equilibrium}}
\date{}
\maketitle
\vspace{-1.5cm}
\begin{center}
{\large 
DONGCHEN LI \\[5pt]
{\it Department of Mathematics, Imperial College London\\
180 Queen's Gate, London SW7 2AZ, United Kingdom}
}
\end{center}
\vspace{1cm}
\footnote{\noindent Date: May 6, 2016.\\
\noindent This work was supported by grant RSF 14-41-00044 at Lobachevsky University of Nizhny Novgorod. The author also acknowledges support by the Royal Society grant IE141468 and EU Marie-Curie IRSES Brazilian-European partnership in Dynamical Systems (FP7-PEOPLE-2012-IRSES 318999 BREUDS).
}
\par{}
\noindent{\bf Abstract.} We show that heterodimensional cycles can be born at the bifurcations of a pair of homoclinic loops to a saddle-focus equilibrium for flows in dimension 4 and higher. In addition to the classical heterodimensional connection between two periodic orbits, we found, in intermediate steps, two new types of heterodimensional connections: one is a heteroclinic between a homoclinic loop and a periodic orbit with a 2-dimensional unstable manifold, and the other connects a saddle-focus equilibrium to a periodic orbit with a 3-dimensional unstable manifold.
\par{}
\par{}
\section{Introduction}
For multidimensional systems, i.e. 4-dimensional flows and 3-dimensional maps, the existence of heterodimensional cycles is a main mechanism that leads to non-hyperbolicity (see \cite{bdv}). A heterodimensional cycle is created by two heteroclinic connections between two saddle periodic orbits with different indices (dimensions of the unstable invariant manifolds). These cycles can be persistent: even when removed by a small perturbation of the system they can re-emerge after an additional arbitrarily small perturbation (see \cite{d92,d95,d95_2,bd96}). Thus, they give a mechanism for a persistent coexistence of saddles with different dimension of the unstable manifold within the same chaotic attractor. These attractors should exhibit properties different from those predicted by hyperbolic theory, e.g. the shadowing property could be violated (see \cite{dgsy27}). Therefore, the study of heterodimensional cycles is important for a further advancement of the theory of multidimensional chaos. 
\par{}
In this paper, we consider multidimensional flows with two Shilnikov loops (i.e. homoclinic loops to a saddle-focus equilibrium) and show under which conditions their bifurcations can create heterodimensional cycles. As intermediate steps, we find two new types of heterodimentional connections. One type of the connection is between a homoclinic loop and a periodic orbit of index 2, and the other one is between a saddle-focus equilibrium and a periodic orbit of index 3; studying these bifurcations could be of independent interest.
\par{}
We consider a $C^r$-flow $X$ in $\mathbb{R}^n$ (where $r\geqslant 3, n \geqslant 4$), which has an equilibrium $O$ of saddle-focus type. 
The eigenvalues of the linearized matrix of $X$ at $O$ are 
\[
\gamma,-\lambda+\omega i,-\lambda-\omega i,\alpha_j 
\]
where Re$(\alpha_j)<-\lambda<0<\gamma$ ($j=1,2\ldots n-3$), $\omega\neq 0$, and we assume
\begin{equation}\label{eq:setting_rho}
\rho=\dfrac{\lambda}{\gamma}<\dfrac{1}{2}.
\end{equation} 
From now on we let $\gamma=1$ (this can always be achieved by time scaling).
It follows from the result in Appendix A of \citep{sstc1} that the system near $O$ can be brought to the form
\begin{equation}\label{eq:setting_1}
\begin{array}{rcl}
\dot{y}&=&y, \\[5pt]
\dot{x_1}&=&-\rho x_1 - \omega x_2 + f_{11}(x,y,z)x+f_{12}(x,y,z)z, \\[5pt]
\dot{x_2}&=&\omega x_1 - \rho x_2 + f_{21}(x,y,z)x+f_{22}(x,y,z)z, \\[5pt]
\dot{z}  &=&Bz+f_{31}(x,y,z)x+f_{32}(x,y,z)z,    
\end{array} 
\end{equation}
\noindent by some $C^{r-1}$-transformations of coordinates and time, where $x=(x_1,x_2)$, and the eigenvalues of matrix $B$ are $\alpha_1\ldots \alpha_{n-3}$.
Functions $f_{ij}$ are $C^{r-1}$-smooth and satisfy
\begin{equation}\label{eq:setting_2}
f_{ij}(0,0,0)=0,\quad f_{1j}(x,0,z)\equiv 0,\quad f_{2j}(x,0,z)\equiv 0, \quad f_{i1}(0,y,0) \equiv 0 \quad(i=1,2,3\mbox{ and }j=1,2).
\end{equation} 
In such coordinate system, the coordinates of $O$ are $(0,0,0)$ and the local invariant manifolds are straightened, i.e., 
$W^u_{loc}(O)=\{x=0,z=0\}, W^s_{loc}(O)=\{y=0\}$ and $W^{ss}_{loc}(O)=\{x=0,y=0\}$.

\begin{figure}[!h]
\begin{center}
\includegraphics[scale=0.45]{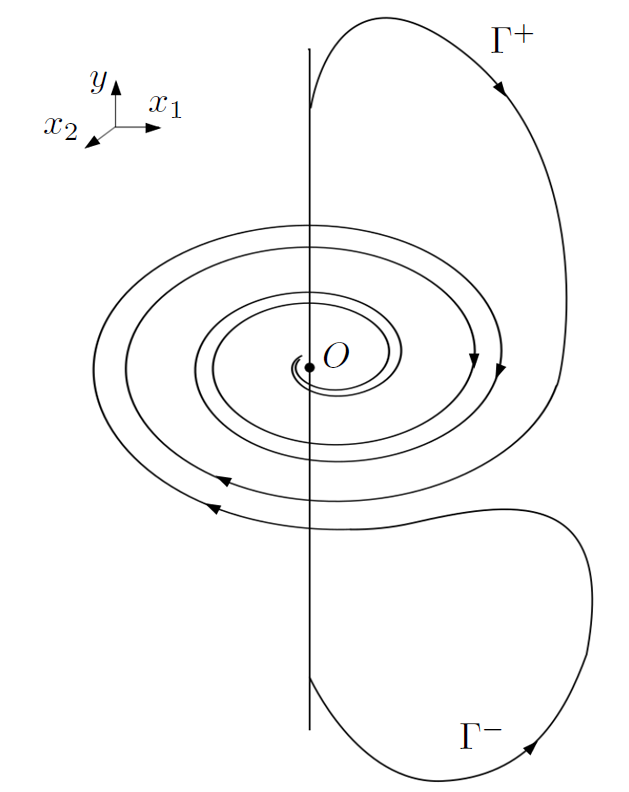}
\end{center}
\caption{The system $X$ has two homoclinic loops $\Gamma^+$ and $\Gamma^-$ to the saddle-focus equilibrium $O$.}
\label{fig:setting}
\end{figure}

The one-dimensional unstable manifold of $O$ consists of two separatrices; the upper one, $\Gamma^+$ corresponds, locally, to $y>0$
and the lower separatrix $\Gamma^-$ corresponds to $y<0$. Let the separatrices $\Gamma^+$ and $\Gamma^-$ return to $O$ as $t\to+\infty$ and form homoclinic
loops. Thus, each of the homoclinic loops, when it tends to $O$ as $t=-\infty$, coincides with a piece of the $y$-axis, and when the loop tends to $O$ as $t\to+\infty$
it lies in the space $y=0$ (see figure \ref{fig:setting}). We assume that $\Gamma^+$ and $\Gamma^-$ do not lie in the strong-stable manifold $W^{ss}$, i.e. $x\neq 0$ as $t\to+\infty$. 
\par{}
Note that two homoclinic loops are necessary for our construction. It was shown in the work (\citep{os87,os92}) of Ovsyannikov and Shilnikov that there exist periodic orbits of different indices near a Shilnikov loop. However, those orbits cannot belong to the same chain-transitive set near the loop. Therefore, no heterodimensional cycle can be born at the bifurcation of one loop. Indeed, assume that the system $X$ has only one homoclinic loop $\Gamma^+$. It is known (see \citep{tu96}) that, under some genericity assumptions on the loop, the system $X$, and every system close to it, have a 3-dimensional invariant manifold $\mathcal{M}$ such that every orbit that lies entirely in the small neighborhood $U$ of $O\cup \Gamma^+$ must lie in $\mathcal{M}$. If there exists a heterodimensional cycle related to two periodic orbits $P$ and $Q$ near the loop $\Gamma^+$, then this cycle must lie in $\mathcal{M}$. This is impossible because, in the 3-dimensional flow on $\mathcal{M}$, the orbit with larger index, say $Q$, becomes completely unstable, which means that there is no heteroclinic connection in $\mathcal{M}$ from $P$ to $Q$. If we want to create a heterodimensional cycle, there should be two saddle periodic orbits with different indices which are not contained in the same 3-dimensional invariant manifold.
\par{}
This situation becomes possible when we consider the bifurcation of two homoclinic loops $\Gamma^+$ and $\Gamma^-$. Even in this case, there can still exist a 3-dimensional invariant manifold containing $\Gamma^+$ and $\Gamma^-$ (see \citep{sstc1}). To avoid this, we need to break the necessary and sufficient condition (proposed in \citep{tu96}) for the existence of a normally hyperbolic invariant manifold that contains both homoclinic loops (see also \citep{bc15}).
\par{}
In order to do this, let us first impose a certain non-degeneracy condition on the system $X$ (this condition is open and dense in $C^r$, i.e., 
if it is not fulfilled initially, then it can be achieved by an arbitrarily small perturbation of the system; once this condition is satisfied, it holds for every $C^r$-close system). Consider an extended unstable manifold $W^{uE}(O)$ of $O$. This is is a smooth 3-dimensional invariant manifold which contains $W^u(O)$ 
and is transverse to $W^{ss}_{loc}(O)$ at $O$. In the coordinates where the system assumes form (\ref{eq:setting_1}), the manifold $W^{uE}_{loc}$ 
is tangent to $z=0$ at the points of $W^u_{loc}$ (see chapter 13 of \citep{sstc2}).\\
\\
{\bf Non-degeneracy Condition}: The extended unstable manifold $W^{uE}(O)$ 
is transverse to the strong-stable foliation of the stable manifold $W^s(O)$ at the points of the homoclinic loops $\Gamma^+$ and $\Gamma^-$.\\
\\
The strong-stable foliation $\mathcal{F}_0$ is the uniquely defined, smooth, invariant foliation of the stable manifold, which includes $W^{ss}(O)$ as one of its leaves;
in the coordinates of (\ref{eq:setting_1}), the leaves of the foliation in a neighborhood of $O$ are given by $(y=0,x=const)$. The transversality condition
implies that the closed invariant set $O\cup\Gamma^+\cup\Gamma^-$ is partially-hyperbolic: at the points of this set the contraction along the 
strong-stable leaves is stronger than a possible contraction in the directions tangent to $W^{uE}$. The partial hyperbolicity implies that the strong-stable foliation $\mathcal{F}_0$ extends, as a locally invariant, absolutely continuous foliation with smooth leaves, to a neighborhood $U$
of $O\cup\Gamma^+\cup\Gamma^-$, and the foliation persists for all $C^r$-close systems (see \citep{an67,ts98}). 
\par{}
Let us take a small cross-section $\Pi$ to the local stable manifold $W^s_{loc}(O)$ such that both $\Gamma^+$ and $\Gamma^-$ intersect $\Pi$. The intersections of the orbits of the leaves of $\mathcal{F}_0$ by the flow with the cross-section $\Pi$ form a strong-stable invariant foliation $\mathcal{F}_1$ for the Poincaré map $T$, which has leaves of the form $(x,y)=h(z)$ where the derivative $h'(z)$ is uniformly bounded. This foliation is invariant in the sense that $T^{-1}(l\cap(T(\Pi)))$ is a leaf of the foliation if the intersection is non-empty. The foliation is also contracting in the sense that, for any two points in the same leaf, the distance of their iterates under the map $T$ tends to zero exponentially. Besides, this foliation is absolutely continuous such that the projection along the leaves from one transversal to another one changes areas by a finite multiple bounded away from zero. One can see \citep{an67} for more discussions on the properties of such foliations. The detailed sufficient condition for the existence of the strong-stable foliation with above-mentioned properties is proposed in \citep{ts98} and our system $X$ satisfies this condition. Note that condition $\rho<\dfrac{1}{2}$
implies that the flow near $O$ expands three-dimensional volume in the $(x,y)$-space; the partial hyperbolicity and the fact that the
orbits in $U$ spend only a finite time between consequent returns to the small neighborhood of $O$ imply that the flow in $U$ uniformly
expands the three-dimensional volume transverse to the strong-stable foliation. Correspondingly, the Poincaré map $T$ expands the two-dimensional
area transverse to the strong-stable foliation on $\Pi$. 
\par
\begin{figure}[!h]
\begin{center}
\includegraphics[scale=0.33]{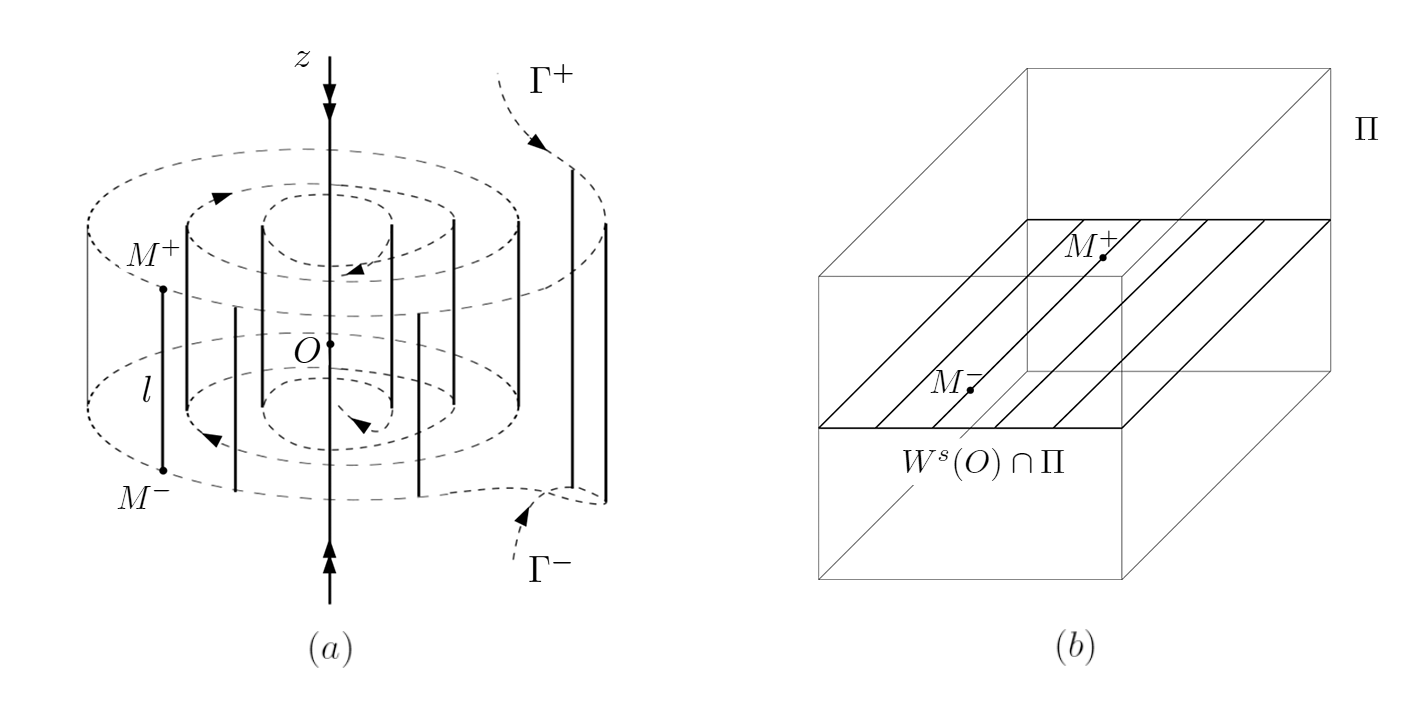}
\end{center}
\caption{The foliation $\mathcal{F}_0$ on the local stable manifold $W^s_{loc}(O)$ is shown schematically in (a), where the dashed curves represent the two homoclinic loops and the solid vertical lines represent leaves of $\mathcal{F}_0$. The coincidence condition for the system $X$ is that, for any point in $\Gamma^+\cap W^s(O)$ lying on a leaf $l$, there exists one point in $\Gamma^-$ that also lies on $l$. For a 4-dimensional flow, the foliation $\mathcal{F}_1$ on a small 3-dimensional cross-section $\Pi$ is shown in (b), where the intersection points $M^+$ and $M^-$ belong to the same leaf.}
\label{fig:coincidence}
\end{figure}

\begin{figure}[!h]
\begin{center}
\includegraphics[scale=0.66]{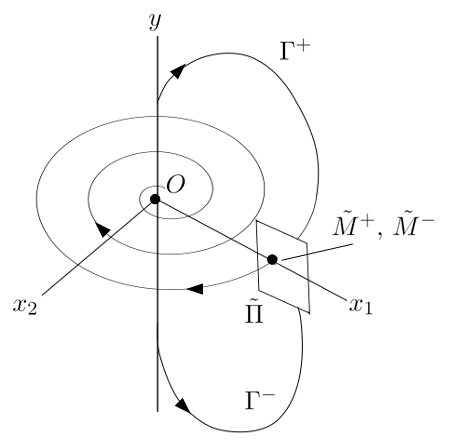}
\end{center}
\caption{For the quotient system that satisfies the coincidence condition, the homoclinic loops $\Gamma^+$ and $\Gamma^-$ intersect the quotient cross-section $\tilde{\Pi}$ at the same point, where  $\tilde{\Pi}$ is $\Pi\cap\{z=z^*\}$ for some $\|z^*\|<\delta$, and $\tilde{M}^+$ and $\tilde{M}^-$ are projections of $M^+$ and $M^-$ on $\tilde{\Pi}$ along the leaves of the strong-stable foliation $\mathcal{F}_1$.}
\label{fig:quotient}
\end{figure}

{}
We can now introduce an assumption on our system $X$ which prevents the 3-dimensional reduction.\\
\\
{\bf Coincidence Condition}: The two separatrices $\Gamma^+$ and $\Gamma^-$ are included in the same set of leaves of the strong-stable foliation $\mathcal{F}_0$ on $W^s(O)$, which means that, for any point $M^+\in\Gamma^+$ lying on a leaf $l$, there exists a corresponding point $M^-\in\Gamma^-$ also lying on the leaf $l$ (see figure \ref{fig:coincidence}). \\\\
This condition means that the projections of $W^s(O)\cap\Gamma^+$ and $W^s(O)\cap\Gamma^-$ onto any transversal along those leaves coincide. Now consider the foliation $\mathcal{F}_1$ on a small cross-section $\Pi$ defined above. Note that leaves of $\mathcal{F}_1$ are obtained by following the orbits of leaves of $\mathcal{F}_0$. Therefore, the intersection points of $\Gamma^+$ and $\Gamma^-$ with $\Pi$ have the same $x$-coordinates (since $\Pi$ is near $O$ and the foliations on $W^s_{loc}(O)$ are straightened). The system satisfying the coincidence condition after taking quotient along the leaves strong-stable foliation $\mathcal{F}_0$ is shown in figure \ref{fig:quotient}. 
\par{}
In this paper, we show how heterodimensional cycles can be created when a system $X$ that satisfies the non-degeneracy and coincidence conditions is perturbed. Note that what we study here is a codimension-3 bifurcation (the existence of two homoclinic loops and the coincidence condition give 3 equality-type conditions imposed on the system). The problem becomes codimension-1 if one considers a class of symmetric systems (then the existence of
one loop implies the existence of the other, and the coincidence condition can be satisfied automatically). In \citep{lt15}, we showed the birth of heterodimensional cycles at such bifurcation, where the system $X$ is invariant with respect to the transformation $\mathcal{R}:(x,y,z)\to (x,-y,\mathcal{S}z)$ where $\mathcal{S}$ is an involution which changes signs of some of the $z$-coordinates. We considered perturbations keeping the symmetry which made the construction quite complicated. As we show in the present paper, different and much simpler constructions become possible when we do not restrict ourselves by symmetric perturbations only. 
\section{Results}
Before stating our main result, let us introduce the parameters. We denote by $M^+$ and $M^-$ the first intersection points of $\Gamma^+$ and $\Gamma^-$ with the cross-section $\Pi$. Let $\zeta$ be a parameter describing the relative position of $M^+$ and $M^-$ in coordinates $x$. Note that, at $\zeta=0$, the system $X$ satisfies the coincidence condition. When we perturb the system, $\zeta$ can become non-zero (i.e. the coincidence condition is no longer fulfilled). Another parameter we use is $\rho$ defined by (\ref{eq:setting_rho}). Below we will consider a 2-parameter family $\{X_{\zeta,\rho}\}$ of flows. We remark here that this 2-parameter family $\{X_{\zeta,\rho}\}$ does not unfold the two homoclinic loops, but it can be viewed as a special choice of parameter values within a 4-parameter unfolding, where two more parameters are used to control the splitting of the loops. 
\begin{thm}\label{thm1}
Let $X_{\zeta,\rho}$ be a family of $C^r$ flows in $\mathbb{R}^n$ ($r\geqslant 3,n \geqslant 4$), where $X_{0,\rho^*}=X$. There exist a sequence $\{\zeta_{i},\rho_i\}$ of parameter values, where $\rho_i\to\rho^*$ and $\zeta_i \to 0$ as $i\to +\infty$, such that each pair $(\zeta_i,\rho_i)$ of parameter values corresponds to a system $X_{\zeta_{i},\rho_i}$ having a heterodimensional cycle with two saddle periodic orbits of indices 2 and 3.
\end{thm}
\par{}
We prove this theorem in the next section. In what follows, we explain the method used in the proof of Theorem \ref{thm1}. To obtain a heterodimensional cycle in $X$, we only need to create one for the Poincaré map $T$ on a cross-section $\Pi$, where this cycle is related to two saddle periodic points $P,Q\in\Pi$ of indices 1 and 2, respectively. The candidates for the index-1 point $P$ are provided by Shilnikov theorem (see \citep{sh65,sh70}), which states that the homoclinic loop $\Gamma^\pm$ is contained in the closure of a hyperbolic set and the intersection of this set with the cross-section $\Pi$ has infinitely many index-1 saddle fixed points $P_k^\pm$ accumulating onto the intersection point $M^\pm$ of $\Gamma^\pm$ with $\Pi$ (see Lemma \ref{lem:WsP_1}). The next step to show is that there exists a saddle periodic point $Q$ of index 2 which has two heteroclinic connections with one of the points $P_k^\pm$. For certainty, we consider the case where we pick $P$ from $\{P_k^+\}$. The same results can be achieved when $P\in\{P_k^-\}$.
\par
\begin{figure}[!h]
\begin{center}
\includegraphics[scale=0.42]{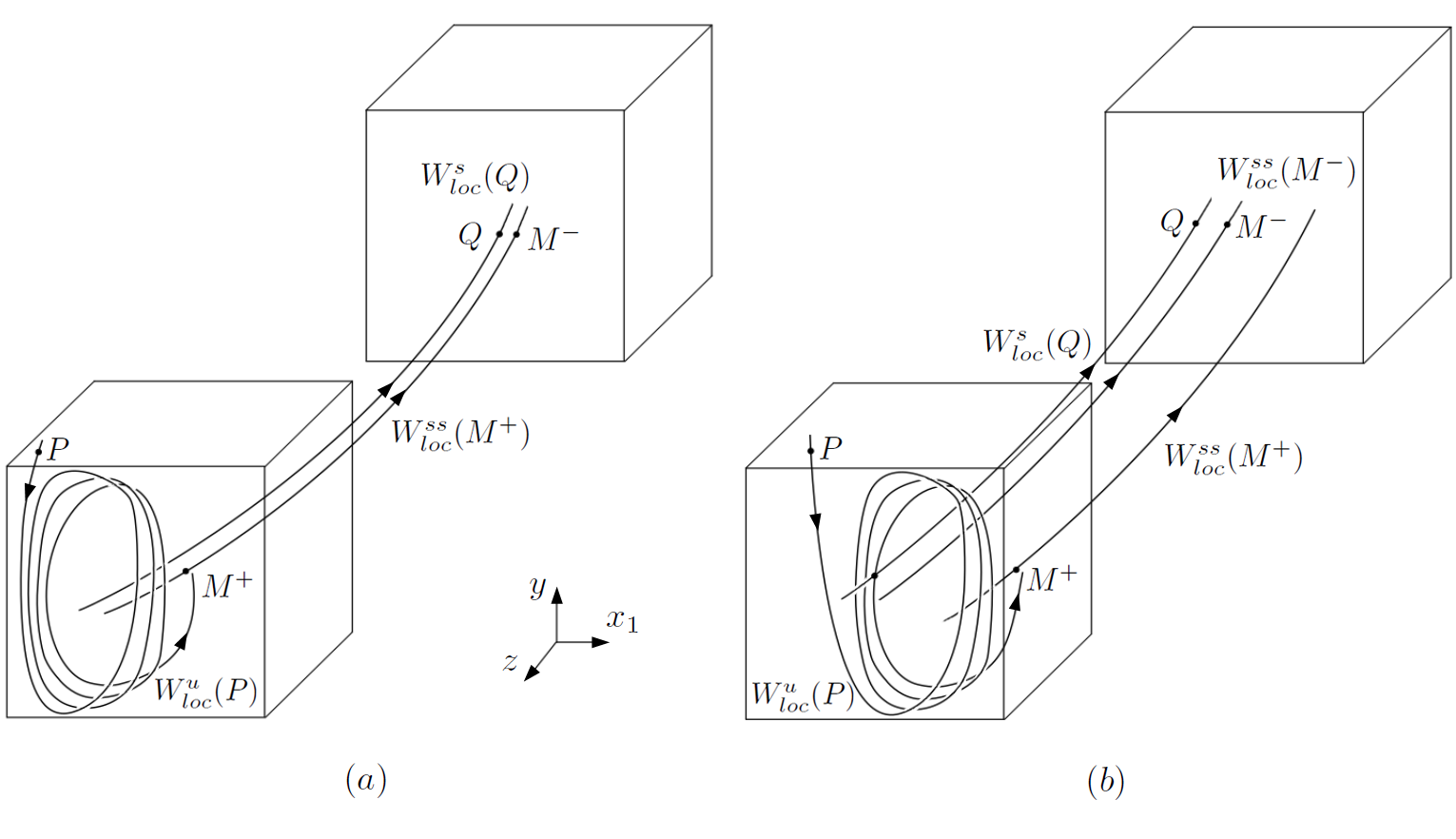}
\end{center}
\caption{The two cubes represent two connected components of the image of cross-section $\Pi$ under the Poincaré map (see (\ref{eq:setting_7}) and (\ref{eq:setting_8})), where $\Pi$ is taken as a small piece of $\{x_2=0\}$. The partial hyperbolicity given by the non-degeneracy condition ensures that the strong-stable foliation on $W^s(O)$ can be extended to a neighborhood of $O\cup\Gamma^+\cap\Gamma^-$, and the strong stable manifolds are leaves of this foliation. As shown in figure (a), we can create an index-2 point $Q$ arbitrarily close to $M^-$ by changing $\rho$ at $\zeta=0$ (so $M^+$ and $M^-$ lie on the same leaf). In figure (b), the intersection $W^s(Q)\cap W^u(P)$ is created by making $\zeta$ non-zero (while changing $\rho$ together to keep $Q$ index-2).}
\label{fig:3d_1}
\end{figure}

{}
We use a modification of the result in \citep{os87} that if a system has a homoclinic loop to a saddle-focus with $\rho<\dfrac{1}{2}$,
then, by an arbitrarily small perturbation which changes the value of $\rho$ without splitting the loop, one can find an index-3 periodic orbit intersecting the cross-section twice, and this orbit can be chosen such that it is as close as we want to the homoclinic loop. By applying this result to the homoclinic loop $\Gamma^-$ in our system $X$, we can change $\rho$ to obtain a saddle periodic point $Q\in\Pi$ of the Poincaré map $T$ with period 2 and index 2 arbitrarily close to $M^-$ (see figure \ref{fig:3d_1} (a) and Lemma \ref{lem:p2i2}). Moreover, by changing $\rho$ and $\zeta$ together, we can make the quasi-transverse intersection $W^s(Q)\cap W^u(P)$ non-empty at the same time (see figure \ref{fig:3d_1} (b)). Here quasi-transversality means that, for two manifolds $U$ and $V$, we have $\mathcal{T}_x U\cap\mathcal{T}_x V=\{0\}$ for the intersection point $x$ of $U\cap V$, where $\mathcal{T}_xU$ and $\mathcal{T}_xV$ are tangent spaces. The next step is to prove that a transverse intersection $W^s(P)\cap W^u(Q)$ also exists at this moment. We achieve this by considering the facts that the map $T$ expands 2-dimensional areas in $(x,y)$-directions and the points in the set $\{P_k^+\}$ including $P$ are homoclinically related (see \ref{prf:transverse1}).
\par{}
Note that, by a sequence of small perturbation in $\zeta$ and $\rho$, we can create a sequence $Q_i$ of index-2 periodic points such that $W^s(Q)\cap W^u(P)\neq\emptyset$, where each of them corresponds to certain pair $(\zeta_i,\rho_i)$ of parameter values and $Q_i\to M^-$. The stable manifolds of these points are given by the leaves of the strong-stable foliation through these points, so we have $W^s(Q_i)\to W^{ss}(M^-)$ as $i\to +\infty$. It follows that we have $W^u(P)\cap W^{ss}(M^-)$ in the limit, which implies a heterodimensional connection between the homiclinic loop $\Gamma^-$ and a periodic orbit of index 2 corresponding to the point $P\in\Pi$. This is a new type of bifurcation similar to "generalised" or "super" homoclinics of \citep{ca10,ekts89,st97,tu01,ho96}. 
\par
\begin{figure}[!h]
\begin{center}
\includegraphics[scale=0.42]{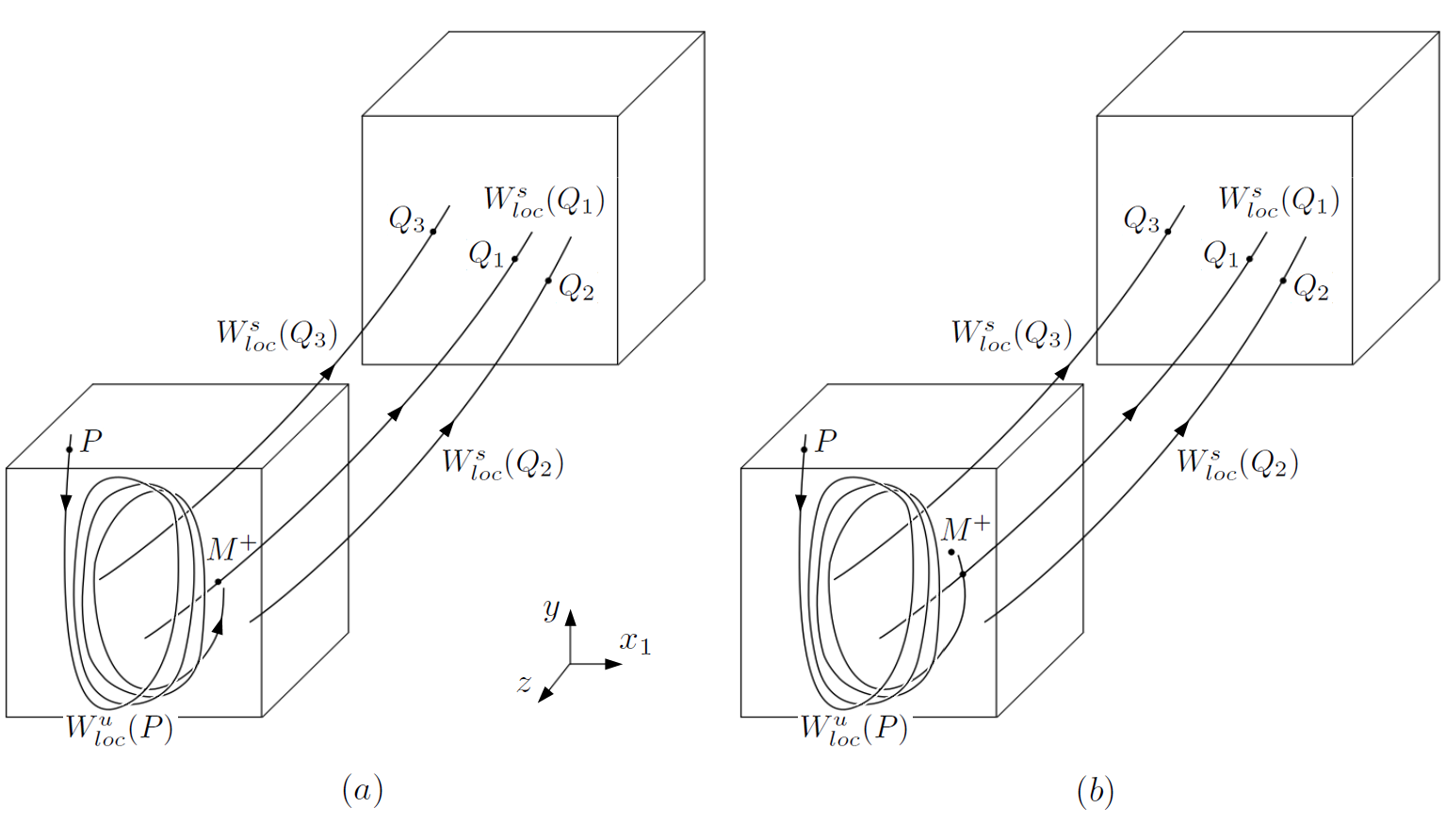}
\end{center}
\caption{Figure (a) shows the critical moment that $M^+\in W^s_{loc}(Q_1)$, where the periodic orbit $\{Q=Q_1,Q_2,Q_3\}$ is of period 3 and index 2. In figure (b), we achieve the intersection $W^s(Q_1)\cap W^u(P)$ by changing $\zeta$.}
\label{fig:3d_2}
\end{figure}

{}
As mentioned in the beginning of this paper, we have another type of bifurcation that creates heterodimensional cycles. In this case, we will split the two homoclinic loops. To do this, we need to introduce two parameters $\mu_1$ and $\mu_2$ that control the splitting of the two homoclinic loops $\Gamma^+$ and $\Gamma^-$, respectively (i.e. the loops split at a non-zero velocity as those parameters change). More specifically, we let $\mu_1$ and $\mu_2$ be the $y$-coordinates of the intersection points $M^+$ and $M^-$, respectively. We still try to create a heterodimensional cycle on $\Pi$ which is related to an index-1 point $P$ and an index-2 point $Q$. Similarly, We choose the index-1 point $P$ from $\{P^+_k\}$ (the same result holds for $P\in\{P^-_k\}$). Then we find a periodic point $Q$ of period 3 and index 2 such that the point $M^+$ falls onto its local stable manifold $W^s_{loc}(Q)$ (see figure \ref{fig:3d_2} (a)). Since the manifold $W^s_{loc}(Q)$ is not straightened, it is not easy to put $M^+$ onto it. In order to do this, we need a freedom to change two more parameters $u$ and $v$ which are smooth functions of coefficients of the system (see (\ref{eq:pp2_26}) and (\ref{eq:pp2_27})). Therefore, what we consider here is a 6-parameter unfolding with $\mu_1,\mu_2,\zeta,\rho,u$ and $v$. Note that the local unstable manifolds $W^u_{loc}(P_k^+)$ of the index-1 fixed points $P_k^+$ given by Shilnikov theorem are spirals winding onto the intersection point $M^+$ (see (\ref{eq:thm2_1})). Hence, by an arbitrarily small perturbation in $\zeta$ (to move $M^+$ out of $W^s_{loc}(Q)$), we can create the quasi-transverse intersection $W^s(Q)\cap W^u(P)$ (see figure \ref{fig:3d_2} (b)). We also show that a transverse intersection $W^s(P)\cap W^u(Q)$ exists at this moment by a similar method used for Theorem \ref{thm1}. This means that a heterodimensional cycle related to $P$ and $Q$ is created.  
\par{}
The result that $M^+\in W^s_{loc}(Q)$ obtained in the intermediate step implies a heterodimensional connection between the saddle-focus equilibrium $O$ and a periodic orbit of index 3 corresponding to the point $Q\in\Pi$. The following result holds.
\begin{spacing}{1.3}\begin{thm}\label{thm2}
We consider a 5-parameter family $X_{\mu,\zeta,\rho,u,v}$ of $C^r$ flows in $\mathbb{R}^n$ ($r\geqslant 3,n \geqslant 4$), where we have $\mu_1=-\mu_2=\mu$ and $X_{0,0,\rho^*,u^*,v^*}=X$. By an arbitrarily small perturbation, we can always make the parameter values $\rho^*,u^*$ and $v^*$ satisfy $\rho^*=\frac{p}{q}\in(\mathbb{Q}\cap(0,\frac{1}{2})),u^*=\frac{p_1}{q}$ and $v^*=\frac{p_2}{q}$, where $p,q$ are co-prime and $p_1,p_2$ are two integers. The triple $(\rho^*,u^*,v^*)$ gives a sequence $\{(\mu_{j},\zeta_j,u_{j}, v_{j})\}$ accumulating on $(0,0,u^*,v^*)$ as $j\to +\infty$ such that the corresponding system $X_{\mu_j,\zeta_j,\rho^*,u_{j}, v_{j}}$ has a heterodimensional cycle related to two saddle periodic orbits with indices 2 and 3. 
\end{thm}\end{spacing}
\par{}
We now give an example of an arbitrarily small perturbation that can be used to obtain the triple $(\rho^*,u^*,v^*)$ stated in the theorem from a general one $(\rho,u,v)$. Let us write $\rho,u$ and $v$ as
\begin{equation}
\begin{array}{rcl}
\rho=a_1a_2\dots a_n\dots\,,\quad
u=b_1b_2\dots b_n\dots\,,\quad
v=c_1c_2\dots c_n\dots\,,
\end{array}
\end{equation}
\noindent and define $\rho^*,u^*$ and $v^*$ as
\begin{equation}
\rho^*=\dfrac{a_1a_2\dots a_N1}{10^{N+1}}\,,\quad
u^*=\dfrac{b_1b_2\dots b_N1}{10^{N+1}}\,,\quad
v^*=\dfrac{c_1c_2\dots c_N1}{10^{N+1}}\,.
\end{equation}
\noindent Obviously, we can get a desired triple $(\rho^*,u^*,v^*)$ by choosing a sufficiently large $N$.
\section{Proofs}
\subsection{Proof of Theorem \ref{thm1}}
The proof is divided into five parts. In part 1, we introduce the family $X_{\zeta,\rho}$ of systems under consideration and describe the Poincaré map $T$ on a cross-section $\Pi$. Then we find periodic points $P$ and $Q$ of the map $T$ with different indices (index 1 and 2) in part 2 and 3, respectively. In the last two parts, we show that, for certain sequence $\{(\zeta_i,\rho_i)\}$ of parameter values, the periodic points $P$ and $Q$ remain index-1 and index-2, and the heteroclinic intersections $W^s(Q)\cap W^u(P)$ and $W^s(P)\cap W^u(Q)$ exist. This gives us a heterodimensional cycle of the map $T$ corresponding to one in the system $X_{\zeta_i,\rho_i}$.  
\subsubsection{Construction of the Poincaré map $T$}
\par{}
Recall that the local stable manifold $W^s_{loc}(O)$ is straightened and has the form $\{y=0\}$. We pick two points $M^+(x^+,0,y^+,z^+)$ and $M^-(x^-,0,y^-,z^-)$ near the equilibrium $O$ such that $M^+\in (\Gamma^+\cap W^s_{loc}(O))$ and $M^-\in (\Gamma^-\cap W^s_{loc}(O))$. We define $\Pi=\{(x_1,0,y,z)\mid|x_1-x^+|<\delta,\|z\|<\delta,|y|<\delta\}$ with upper part $\Pi_1:=\Pi\cap\{ y\geqslant 0\}$ and lower part $\Pi_2:=\Pi\cap\{ y\leqslant 0\}$. Denote by $\Pi_0$ the intersection of $\Pi$ with $W^s_{loc}(O)$. Note that $\|z\|$ decreases much faster than $\|x\|$ along the homoclinic loops as $t\to+\infty$ so that we can assume $\|z^+\|<\delta$. Points $M^+$ and $M^-$ are the intersection points of $\Gamma^+$ and $\Gamma^-$ with $\Pi$, and have coordinates $(x^+,y^+,z^+)$ and $(x^-,y^-,z^-)$ on it. Let $\zeta=x^- -x^+$. We now consider a family $X_{\zeta,\rho}$ of perturbed systems, so $x^\pm,y^\pm$ and $z^\pm$ are smooth functions of parameters. For our system $X$, we have $\zeta=0$ by the coincidence condition.
\par{}
In order to obtain the formula for the Poincaré map $T$, we need the help of two global cross-sections $\Pi_{glo_1}=\{(x_1,x_2,y=d,z)\}$ and $\Pi_{glo_2}=\{(x_1,x_2,y=-d,z)\}$, where $d>0$. The Poincaré map $T$ restricted to $\Pi_i$ $(i=1,2)$ is the composition of a local map $T_{loc_i}:\Pi_i\to \Pi_{glo_i} \, , \,(x_0,y_0,z_0) \mapsto (x_1,x_2,z_1)$ and a global map $T_{glo_i}:\Pi_{glo_i}\to \Pi \, , \,(x_1,x_2,z_1) \mapsto (\bar{x_0},\bar{y_0},\bar{z_0})$. The map $T_{loc_i}$ is given by (see page 738 of \citep{sstc2})
\begin{equation}\label{eq:setting_3}
\begin{array}{l}
 x_{1}=x_{0}\bigg(\dfrac{y_{0}}{(-1)^{i+1}d}\bigg)^{\rho }\cos \left( \omega \ln \left( \dfrac{(-1)^{i+1}d}{y_{0}}
\right) \right)+o\left( |y_{0}|^{\rho } \right) \,,\\[20pt]
 x_{2}=x_{0}\bigg(\dfrac{y_{0}}{(-1)^{i+1}d}\bigg)^{\rho }\sin \left( \omega \ln \left( \dfrac{(-1)^{i+1}d}{y_{0}}
\right) \right)+o(|y_{0}|^{\rho })\,, \\[20pt]
z_1=\left(
\begin{array}{c}
o(|y_{0}|^{\rho })\\
\cdots\\
o(|y_{0}|^{\rho })
\end{array} \right)\,,
\end{array}
\end{equation}
\noindent where the small terms $o(|y|^\rho)$ (for both $y>0$ and $y<0$) are functions of $x,y,z,\varepsilon=(\zeta,\rho)$ satisfying
\begin{equation}\label{eq:smallterms}
\dfrac{\partial^{i+j+k+l} o(|y|^\rho)}{\partial^i y \,\partial^j x \,\partial^k z\,\partial^l \varepsilon}=o(|y|^{\rho-i}) \quad i+j+k+l\leqslant(r-1)\,.
\end{equation}
The global maps $T_{glo_1}$ and $T_{glo_2}$ are diffeomorphisms and can be written in Taylor expansions. We have  
\begin{equation}\label{eq:setting_4}
T_{glo_1}:\quad\left\{
\begin{array}{l}
 \bar{x}_{0}=x^+ +a_{11}x_{1}+a_{12}x_{2}+a_{13}z_1+o\left( |x_{1},x_{2},z_1| \right)\,, \\[8pt]
 \bar{y}_{0}=a_{21}x_{1}+a_{22}x_{2}+a_{23}z_1+o\left(
|x_{1},x_{2},z_{1}| \right) \,, \\[8pt]
 \bar{z}_{0}=z^+ +

 \begin{pmatrix}
 a_{31}x_{1}+a_{32}x_{2}+a_{33}z_1 \\
 \cdots \\
 a_{n-1,1}x_{1}+a_{n-1,2}x_{2}+a_{n-1,3}z_1
 \end{pmatrix}
 
+o\left( |x_{1},x_{2},z_1| \right)\,,
 \end{array}
\right.
\end{equation}
\noindent and
\begin{equation}\label{eq:setting_5}
T_{glo_2}:\quad\left\{
\begin{array}{l}
 \bar{x}_{0}=x^+ +\zeta +b_{11}x_{1}+b_{12}x_{2}+b_{13}z_1+o\left( |x_{1},x_{2},z_1| \right)\,, \\[8pt]
 \bar{y}_{0}=b_{21}x_{1}+b_{22}x_{2}+b_{23}z_1+o\left(
|x_{1},x_{2},z_{1}| \right) \,, \\[8pt]
 \bar{z}_{0}=z^- +

 \begin{pmatrix}
 b_{31}x_{1}+b_{32}x_{2}+b_{33}z_1 \\
 \cdots \\
 b_{n-1,1}x_{1}+b_{n-1,2}x_{2}+b_{n-1,3}z_1
 \end{pmatrix}
 
+o\left( |x_{1},x_{2},z_1| \right)\,,
 \end{array}
 \right.
\end{equation}
\noindent where $a_{j3}$ and $b_{j3}$ $(j=1\cdots n-1)$ are $(n-3)$-dimensional vectors. The Poincaré map $T:\Pi \to \Pi$ is continuous except on points with $y=0.$ Let $T_1:=T|_{\Pi_1}=T_{glo_1}\circ T_{loc_1}:\Pi_1 \to \Pi$ and $T_2:=T|_{\Pi_2}=T_{glo_2}\circ T_{loc_2}:\Pi_2 \to \Pi$. Note that \\[-10pt]
\begin{equation}\label{eq:setting_6}
\lim_{M \to \{y=0\}^+}T_1(M)=M^+ \quad  \mbox{and} \quad \lim_{M \to \{y=0\}^-}T_2(M)=M^-.
\end{equation}
\noindent By the scaling $x=\dfrac{x_0}{x^+}$ and $y=\dfrac{y_0}{d}$, and replacing $\dfrac{\zeta}{x^+}$ by $\zeta$, \noindent the maps $T_1$ and $T_2$ take the form \\[7pt]
\begin{equation}\label{eq:setting_7}
T_1:\quad\left\{\begin{array}{l}

                \bar{y}=Axy^\rho\cos{(\omega \ln {\dfrac{1}{y}}+\eta)}+o(y^\rho) \\[7pt]
                \bar{x}=1+A_1xy^\rho\cos{(\omega \ln {\dfrac{1}{y}}+\eta_1)}+o(y^\rho) \\[7pt]
                \bar{z}=z^+ +
                \begin{pmatrix}
                 A_2xy^\rho\cos{(\omega \ln {\dfrac{1}{y}}+\eta_2)}+o(y^\rho) \\
                 \cdots \\
                  A_{n-2}xy^\rho\cos(\omega \ln {\dfrac{1}{y}}+\eta_{n-2})+o(y^\rho) 
                \end{pmatrix}                
                \end{array}\right.
\end{equation}
\noindent and
\begin{equation}\label{eq:setting_8}
T_2:\quad\left\{\begin{array}{l}
                \bar{y}=-Bx|y|^\rho\cos{(\omega \ln {\dfrac{1}{|y|}}+\theta)}+o(|y|^\rho) \\[7pt]
                \bar{x}=1+\zeta+B_1x|y|^\rho\cos{(\omega \ln {\dfrac{1}{|y|}}+\theta_1)}+o(|y|^\rho) \\[7pt]
                \bar{z}=z^- + 
                \begin{pmatrix}
                 B_2x|y|^\rho\cos{(\omega \ln {\dfrac{1}{|y|}}+\theta_2)}+o(|y|^\rho) \\
                 \cdots \\
                  B_{n-2}x|y|^\rho\cos(\omega \ln {\dfrac{1}{|y|}}+\theta_{n-2})+o(|y|^\rho) 
                \end{pmatrix}     
                 \end{array}\right.,
\end{equation}\\
\noindent respectively, where $z \in \mathbb{R}^{n-3}$, $A=x^+\sqrt{a_{21}^{2}+a_{22}^{2}}$, $
A_{1}=\sqrt{a_{11}^{2}+a_{12}^{2}}$, $ A_{m}=x^+\sqrt{a_{m+1,1}^{2}+a_{m+1,2}^{2}}$ $(m=2,\cdots,n-2)$, $B=x^+\sqrt{b_{21}^{2}+b_{22}^{2}}$, $
B_{1}=\sqrt{b_{11}^{2}+b_{12}^{2}}$, $ B_{m}=x^+\sqrt{b_{m+1,1}^{2}+b_{m+1,2}^{2}}$, $\tan\eta$ $
=-\dfrac{a_{22}}{a21}$, $\tan\eta_{1}=-\dfrac{a_{12}}{a11}$,  $\tan\eta_{m}=-\dfrac{a_{m+1,2}}{a_{m+1,1}}$, ${\tan}\theta$ $
=-\dfrac{b_{22}}{b21}$, $\tan\theta_{1}=-\dfrac{b_{12}}{b11}$,  $\tan\theta_{m}=-\dfrac{b_{m+1,2}}{b_{m+1,1}}$, and the small terms $o(|y|^\rho)$ satisfy (\ref{eq:smallterms}). 
\par{}
From now on, we will work with the maps $T_1$ and $T_2$. Note that the above-mentioned non-degeneracy condition is equivalent to
\begin{equation}\label{eq:nondege}
AA_1\sin(\eta_1-\eta)\neq 0 \mbox{ and } BB_1\sin(\theta_1-\theta)\neq 0.
\end{equation}
\noindent Indeed, in the coordinate system satisfying (\ref{eq:setting_1}) and (\ref{eq:setting_2}), the transversality stated in the non-degeneracy condition is equivalent to the transversality of $T_{glo_1}(\Pi_{glo_1}\cap W^{uE}_{loc}(O))$ and $T_{glo_2}(\Pi_{glo_2}\cap W^{uE}_{loc}(O))$ to the leaves $\{y=0,x=x^+\}$ through $M^+$ and $\{y=0,x=x^-\}$ through $M^-$, respectively, where the extended unstable manifold $W^{uE}_{loc}(O)$ is an invariant manifold tangent to the $\{z=0\}$ (see \citep{sstc1}).  By the formulas (\ref{eq:setting_4}) and (\ref{eq:setting_5}), this is
\begin{equation*}
\mbox{det}\dfrac{\partial(\bar{x}_0,\bar{y}_0)}{(x_1,x_2)}\neq 0
\end{equation*} 
\noindent for both maps $T_{glo_1}$ and $T_{glo_2}$, which is equivalent to 
\begin{equation*}
\begin{array}{rcl}
\begin{vmatrix}
a_{11} & a_{12} \\
a_{21} & a_{22} \\
\end{vmatrix}=AA_1\sin(\eta_1-\eta)\neq 0
& \mbox{ and } &
\begin{vmatrix}
b_{11} & b_{12} \\
b_{21} & b_{22} \\
\end{vmatrix}=BB_1\sin(\theta_1-\theta)\neq 0
\end{array}\,.
\end{equation*}
\subsubsection{Existence of the index-1 fixed point $P$}
Throughout the rest of this paper, we write the coordinates $x,y,z$ in the order $y,x,z$ for being consistent with the way we write the maps $T_1$ and $T_2$.
\par{}
As mentioned before, there exist two countable sets $\{P_k^+\}\subset\Pi_1$ and $\{P_k^-\}\subset\Pi_2$, given by Shilnikov theorem (\citep{sh70}), of index-1 fixed points of $T_1$ and $T_2$ accumulating on $M^+$ and $M^-$, respectively. Indeed, the points $P_k^+$ and $P_k^-$ are obtained by solving the equations $T_1(y,x,z)=(y,x,z)$ and $T_2(y,x,z)=(y,x,z)$, respectively (see below). We now pick an arbitrary point $P$ from the set $\{P_k^+\} \cup \{P_k^-\}$, and we will show that this point $P$ can be the desired index-1 point to create a heterodimensional cycle of $T$. For certainty, we fix a point $P$ from $\{P_k^+\}$. Similar results hold for $P\in \{P_k^-\}$.
\par{}
Let us now find the points $P^+_k$ and their local stable manifolds $W^s_{loc}(P^+_k)$ which will be used later to create the transverse intersection $W^s(P)\cap W^u(Q)$, where $Q$ is an index-2 periodic point. We have the following result.
\begin{lem}\label{lem:WsP_1}
The local stable manifolds $W^s_{loc}(P_k^+)$ of the index-1 fixed points $P_k^+$ given by Shilnikov theorem are graphs of functions $g(x,z)$ defined for all $x$ and $z$ values in $\Pi$, which accumulate, in $C^0$-topology, on $\Pi_0=\{y=0\}$ as $k \to +\infty$.
\end{lem}
\par{}
\noindent \textit{Proof.} We first find the fixed points $P^+_k$, which can be done by plugging $(\bar{y}=y,\bar{x}=x,\bar{z}=z)$ into (\ref{eq:setting_7}). From the last two equations in (\ref{eq:setting_7}), the coordinates $x$ and $z$ can be expressed by $y$, which gives a equation for coordinate $y$:
\begin{equation}\label{eq:WsP_1}
y=Ay^\rho\cos(\omega\ln(\dfrac{1}{y})+\eta)+o(y^\rho).
\end{equation}
We have the fixed points $P^+_k$ with
\begin{equation}\label{eq:WsP_2}
\begin{array}{rcl}
y_k&=&C\exp\Big(\dfrac{-\pi k}{\omega}\Big)+o\Big(\exp\Big(\dfrac{-\pi k}{\omega}\Big)\Big)\, ,\\[10pt]
x_k&=&1+o\Big(\exp\Big(\dfrac{-\pi k}{\omega}\Big)\Big)\, ,\\[10pt]
z_k&=&z^+ +o\Big(\exp\Big(\dfrac{-\pi k}{\omega}\Big)\Big)\, ,
\end{array}
\end{equation}
\noindent where $y_k,x_k,z_k$ are the coordinates of $P^+_k$, $C=\exp\Big(\dfrac{2\eta-\pi}{2\omega}\Big)$, and $k$ is any positive integer greater than some sufficiently large $K$. Note the the points $P^+_k$ are of index 1. Indeed, by Lemma \ref{lem:index_2} (i.e. lemma 5 in \citep{lt15}), we have that a fixed point $P(y,x,z)$ of $T_i$ $(i=1,2)$ is of index 2 only if $\cos\xi$ is bounded away from 0, where $\xi=\omega\ln\dfrac{1}{y}+\eta-2\pi j=\xi$ (see (\ref{eq:pp1_3}) for details). However, the first equation in (\ref{eq:WsP_2}) implies that $\cos\xi_k$ is small when $k$ is sufficiently large. We also note that, under our consideration, the index of a periodic point is at most 2 since the multipliers corresponding to $z$ coordinates stay inside the unit circle during all the small perturbations. 
\par{}
We now consider the inverse image under $T_1$ of a small piece of the surface $\{y=y_k\}$ containing $P^+_k$. By (\ref{eq:setting_7}), we have
\begin{equation}\label{eq:WsP_3}
\sin\Big(\dfrac{\pi}{2}-\theta-\omega\ln\dfrac{1}{y}\Big)=\dfrac{1}{xA}\Big(\dfrac{y_k}{y^\rho}+o(1)_{y\to 0}\Big)\, ,
\end{equation}
\noindent where $(y,x)$ are coordinates of the points in the inverse image ($z$ coordinates are in the small term) and $x$ is bounded since the small cross-section $\Pi$ is bounded. We have following equation if $y$ and $\dfrac{y_k}{y^\rho}$ are sufficiently small:
\begin{equation}\label{eq:WsP_4}
\dfrac{\pi}{2}-\theta-\omega\ln\dfrac{1}{y}=\dfrac{1}{xA}\Big(\dfrac{y_k}{y^\rho}+o(1)_{y\to 0}\Big)+m\pi\quad ,m=0,\pm1,\pm2,\ldots \quad,
\end{equation}
\noindent which, by noting that the surface contains $P^+_k$, leads to
\begin{equation}\label{eq:WsP_5}
y=C\exp(\dfrac{-\pi k}{\omega})+o(1)_{k\to +\infty}\, .
\end{equation}
\noindent Formula (\ref{eq:WsP_5}) is valid for all values of $x,z$, where $(y,x,z)\in\Pi$, if $y$ and $\dfrac{y_k}{y^\rho}$ are sufficiently small. This requirement is equivalent to that $k$ is sufficiently large. One can check that the successive backward iterates of a small piece of the surface $y=y_k$ containing $P^+_k$ take the form as (\ref{eq:WsP_5}), where the term $o(1)$ stays uniformly small. Since $W^s_{loc}(P)$ is the limit of a sequence of those iterates, $W^s_{loc}(P)$ is given by (\ref{eq:WsP_5}). \qed
\subsubsection{Existence of the index-2 periodic point $Q$}
We will find a periodic point $Q$ of $T_2$ having period 2 and index 2. Let us first introduce a transformation for $y$-coordinates of points on $\Pi_2$: 
\begin{equation}\label{eq:pp1_3}
\omega\ln\dfrac{1}{|y|}=2\pi j+ \xi -\theta \, ,\quad\xi\in[0,2\pi)\,,
\end{equation}
\noindent by which we divide the cross-section into different regions with index $j$ and let $\xi$ be a new coordinate in each region. Thus, given a period-2 orbit $\{Q=Q_1(y_1,x_1,z_1),Q_2(y_2,x_2,z_2)\}$, we define the integers $j_1$ and $j_2$ by the rule
\begin{equation}\label{eq:dense_1}
\omega\ln\dfrac{1}{|y_i|}=2\pi j_i+ \xi_i -\theta \, ,\quad\xi_i\in[0,2\pi)\quad i=1,2\,.
\end{equation}
\noindent We have the following lemma.
\begin{lem}\label{lem:p2i2}
There exists certain function $\Psi(\zeta,\rho,j_1,j_2,c)$ which is uniformly bounded and smooth with respect to $\zeta$ and $c$ such that if the relation
\begin{equation}\label{eq:p2i2}
\rho j_1 - j_2=\Psi(\zeta,\rho,j_1,j_2,c)\,,
\end{equation}
\noindent is satisfied for and $c\in(-1,1)$ and some sufficiently large integers $j_1$ and $j_2$ such that $\dfrac{j_2}{j_1}\in(0,\dfrac{1}{2})$, then the Poincaré map $T$ has an index-2 periodic orbit $\{Q_1(y_1,x_1,z_1),Q_2(y_2,x_2,z_2)\}\subset\Pi_2$ corresponding to $(j_1,j_2)$ via the transformation (\ref{eq:pp1_3}). By taking $j_1$ and $j_2$ sufficiently large, $Q_1$ and $Q_2$ can become arbitrarily close to $M^-$.
\end{lem}
\par{}
Lemma \ref{lem:p2i2} ensures that, by taking $j_1,j_2$ larger and larger, one can create a sequence $\{Q_k\}$ of index-2 periodic points accumulating on $M^-$ (with different parameter values for each $Q_k$).
\par{}
\noindent{\it Proof of Lemma \ref{lem:p2i2}.} We first state a result on the condition for a periodic point to have index 2.
\begin{lem}\label{lem:index_2}
Let a point $Q$ have period k under $T_2$ with the orbit $\{Q=Q_1(y_1,x_1,z_1),$ $Q_2(y_2,x_2,z_2)$ $\cdots,Q_k(y_k,x_k,z_k)\}$. By the transformation above, we have the following relation:
\begin{equation*}
\omega\ln\dfrac{1}{|y_i|}=2\pi j_i+ \xi_i -\theta \quad i=1,2,\cdots,k\,.
\end{equation*} 
The point $Q$ is of index 2, if and only if
\begin{equation*}
\cos(\xi_1-\varphi)\cos(\xi_2-\varphi)\cdots\cos(\xi_k-\varphi)=c\psi(\xi,j,x,z)\,,
\end{equation*}
\noindent where $|c|<1$, $\varphi=\arctan(\dfrac{\omega}{\rho})$, $\xi=(\xi_1,\xi_2,\cdots,\xi_k)$, $j=(j_1,j_2,\cdots,j_k)$, $x=(x_1,x_2,\cdots,x_k)$, $z=(z_1,z_2,\cdots,z_k)$ and $\psi=o(1)_{j_1,j_2,\cdots,j_k \to \infty}$ is certain function depending continuously on $\xi,x,z$ and parameter $\varepsilon$ such that 
\begin{equation*}
\dfrac{\partial^{i+k+l+n} \psi}{\partial^i {\xi} \,\partial^k {x} \,\partial^l {z}\,\partial^n \varepsilon}=o(1)_{j_1,j_2,\cdots,j_k \to \infty} \quad i+k+l+n\leqslant(r-2).
\end{equation*}
\end{lem}
This result is Lemma 5 in \citep{lt15} and we omit the proof here. By Lemma \ref{lem:index_2}, formula (\ref{eq:setting_8}) for $T_2$ and the rule (\ref{eq:dense_1}), an index-2 periodic orbit $\{Q_1,Q_2\}$ is given by the following equations:
{\allowdisplaybreaks
\begin{align}\label{eq:dense_2}
\qquad\qquad  &y_2=-Bx_1|y_1|^\rho\cos{\xi_1}+o(|y_1|^\rho)\,, \\[7pt]
               & x_2=1+\zeta+B_1x_1|y_1|^\rho\cos{(\xi_1+\theta_1-\theta)}+o(|y_1|^\rho) \,,\label{eq:x_2}\\[7pt]
                &z_2=z^- + 
                \begin{pmatrix}
                 B_2x_1|y_1|^\rho\cos{(\xi_1+\theta_2-\theta)}+o(|y_1|^\rho) \\
                 \cdots \\
                  B_{n-2}x_1|y_1|^\rho\cos(\xi_1+\theta_{n-2}-\theta)+o(|y_1|^\rho) 
                \end{pmatrix}    \,,\label{eq:z_2}\\[7pt]
               & y_1=-Bx_2|y_2|^\rho\cos{\xi_2}+o(|y_2|^\rho)\,, \\[7pt]
              &  x_1=1+\zeta+B_1x_2|y_2|^\rho\cos{(\xi_2+\eta_1)}+o(|y_2|^\rho) \,,\label{eq:x_1}\\[7pt]
              &  z_1=z^- + 
                \begin{pmatrix}
                 B_2x_2|y_2|^\rho\cos{(\xi_2+\theta_2-\theta)}+o(|y_2|^\rho) \\
                 \cdots \\
                  B_{n-2}x_2|y_2|^\rho\cos(\xi_2+\theta_{n-2}-\theta)+o(|y_2|^\rho) 
                \end{pmatrix}    \,,\label{eq:z_1}\\[7pt]
               & \cos(\xi_1-\varphi)\cos(\xi_2-\varphi)=c\psi\,.
\end{align}}
\noindent
where $-1<c<1$, $\psi$ is certain function of $\xi_i,j_i,x_i,z_i$ $(i=1,2)$ depending continuously on parameters and $\psi\to 0$ as $j_1,j_2\to +\infty$. We can express $x$ and $z$ as functions of $y$ and get a reduced system given by 
\begin{eqnarray}
&y_2=-B_{\mbox{new}}|y_1|^\rho\cos{\xi_1}+o(|y_1|^\rho)+O(|y_1|^\rho |y_2|^\rho)\,, \label{eq:dense_2.1}\\[7pt]            
               & y_1=-B_{\mbox{new}}|y_2|^\rho\cos{\xi_2}+o(|y_2|^\rho)+O(|y_1|^\rho |y_2|^\rho)\,, \label{eq:dense_2.2}\\[7pt]             
               & \cos(\xi_1-\varphi)\cos(\xi_2-\varphi)=c\psi\,, \label{eq:dense_2.3}
\end{eqnarray}
\noindent where $B_{\mbox{new}}=(1+\zeta)B$ and we drop the subscript for simplicity. Note that it will be shown later in the computation (see (\ref{eq:dense_7})) that $y_1$ and $y_2$ satisfy the relation
\begin{equation*}
|y_1|^\rho\sim y_2\,.
\end{equation*}
\noindent Thus, we replace $o(|y_1|^\rho)+O(|y_1|^\rho |y_2|^\rho)$ and $o(|y_2|^\rho)+O(|y_1|^\rho |y_2|^\rho))$ in equations (\ref{eq:dense_2.1}) and (\ref{eq:dense_2.2}) by $o(|y_1|^\rho)$ and $o(|y_2|^\rho)$, respectively. By applying the rule (\ref{eq:dense_1}) to equations (\ref{eq:dense_2.1}) - (\ref{eq:dense_2.3}), we obtain
\begin{eqnarray}
&\exp\bigg(\dfrac{-2\pi j_2-\xi_2+\theta}{\omega}\bigg)=B\exp\bigg(\dfrac{-2\pi\rho j_1-\rho\xi_1+\rho\theta}{\omega}\bigg)\cos\xi_1+o\bigg(\exp\bigg(\dfrac{-2\pi\rho j_1}{\omega}\bigg)\bigg) \,,\label{eq:dense_3} \\[7pt]
&\exp\bigg(\dfrac{-2\pi j_1-\xi_1+\theta}{\omega}\bigg)=B\exp\bigg(\dfrac{-2\pi\rho j_2-\rho\xi_2+\rho\theta}{\omega}\bigg)\cos\xi_2+o\bigg(\exp\bigg(\dfrac{-2\pi\rho j_2}{\omega}\bigg)\bigg) \,, \label{eq:dense_4}\\[7pt]
& \cos(\xi_1-\varphi)\cos(\xi_2-\varphi)=c\psi\,.\label{eq:dense_5}
\end{eqnarray}
\quad\quad We solve this system with sufficiently large $j_1$ and $j_2$. Note that there are three equations and two variables $\xi_1$ and $\xi_2$, so the solvability of this system will impose a constraint of its parameters, which, as we will show, is equation (\ref{eq:p2i2}). From now on, we denote by dots the small terms which are functions of $\xi_1,\xi_2,j_1,j_2$ and tend to zero as $j_1$ and $j_2$ tend to positive infinity. 
\par{}
Equation (\ref{eq:dense_5}) implies that one of the terms $\cos(\xi_1-\varphi)$ and $\cos(\xi_2-\varphi)$ must be small. Here we assume that  $\cos(\xi_1-\varphi)$ is small and $\cos(\xi_2-\varphi)$ is bounded from zero, from which we obtain
\begin{equation}\label{eq:dense5.5}
\cos(\xi_1-\varphi)=c\psi_1 \,,
\end{equation}
\noindent where $\psi_1$ is certain function of $\xi_1,\xi_2,j_1,j_2$ depending continuously on $\xi_1,\xi_2$ and parameters. Consequently, we get
\begin{equation}\label{eq:dense_9}
\xi_1=\arccos(c\psi_1)+k_1\pi+\varphi=\dfrac{\pi}{2}+k_1\pi+\varphi+\psi_2(\xi_1,\xi_2,j_1,j_2,c) \,,
\end{equation}
\noindent where $\varphi=\arctan(\dfrac{\rho}{\omega})\in(0,\dfrac{\pi}{2})$, $\psi_2=o(1)_{j_1,j_2\to +\infty}\,$ depends continuously on all arguments and parameters, and $k_1=0,1$ since $\xi_1\in[0,2\pi)$. Note that $\psi_2$ varies slightly when we change $c$ from $-1$ to $1$. 
\par{}
We now look at equation (\ref{eq:dense_3}), which gives another expression for $\cos\xi_1$:
\begin{equation}\label{eq:dense_6}
\cos\xi_1=B^{-1}\exp\bigg(\dfrac{2\pi(\rho j_1-j_2)+\theta-\rho\theta+\rho\xi_1-\xi_2}{\omega}\bigg)+\cdots\,.
\end{equation}
\noindent Since equation (\ref{eq:dense_9}) implies that $\cos\xi_1$ is bounded from zero, the sign of $\cos\xi_1$ is the same as that of the first term on the right hand side of (\ref{eq:dense_6}), which is positive. Therefore, we have $k_1=1$ in (\ref{eq:dense_9}). 
\par{}
Let us now find $\xi_2$. From equation (\ref{eq:dense_6}), we have
\begin{equation}\label{eq:dense_7}
\rho j_1-j_2=\omega\ln(B\cos\xi_1)-\theta+\rho\theta-\rho\xi_1+\xi_2+\cdots\,,
\end{equation}
\noindent By noting that $\cos\xi_1$ is finite from (\ref{eq:dense_9}), equation (\ref{eq:dense_7}) implies that $\rho j_1-j_2$ is bounded, so $\rho j_2-j_1$ is large. We divide both sides of (\ref{eq:dense_4}) by $\exp\Big(\dfrac{-2\pi\rho j_2}{\omega}\Big)$ and take the limit $j_1,j_2\to+\infty$. This will give us a solution to (\ref{eq:dense_4}) as
\begin{equation}\label{eq:dense_8}
\xi_2=\dfrac{\pi}{2}+k_2\pi+\cdots \,,
\end{equation}
\noindent where $k_2=0,1$ since $\xi_2\in[0,2\pi)$. It follows that $\cos(\xi_2-\varphi)$ is bounded away from zero, which agrees with our assumption used to obtain (\ref{eq:dense5.5}) from (\ref{eq:dense_5}).
\par{}
Note that, by implicit function theorem, we can express $\xi_1$ and $\xi_2$ as functions of $j_1,j_2,\zeta$ and $\rho$ from (\ref{eq:dense_9}) and (\ref{eq:dense_8}). By plugging the new expressions of $\xi_1$ and $\xi_2$ into (\ref{eq:dense_7}), we obtain the relation (\ref{eq:p2i2}):
\begin{equation*}
\rho j_1 - j_2=\Psi(\zeta,\rho,j_1,j_2,c),
\end{equation*}
\noindent where $\Psi$ is continuous in $\zeta,\rho$ and $c$, and is uniformly bounded. Note that equation \eqref{eq:dense_2.3} given by Lemma \ref{lem:index_2} requires $\rho$ to be in $(0,1/2)$. Also  the relation (\ref{eq:p2i2}) can be rewritten as 
\begin{equation*}
\rho =\dfrac{j_2}{j_1}+\dfrac{\Psi(\zeta,\rho,j_1,j_2,c)}{j_1},
\end{equation*}
\noindent which implies
\begin{equation}
\rho\sim\dfrac{j_2}{j_1}.
\end{equation}
Therefore, we need to consider $j_1$ and $j_2$ which satisfy not only (\ref{eq:p2i2}) but also $j_2/j_1\in(0,1/2)$. Each such pair $(j_1,j_2)$ gives an index-2 periodic orbit $\{Q_1,Q_2\}$ of $T_2$.
\par{}
By taking $j_1$ and $j_2$ sufficiently large with $j_2/j_1\in(0,1/2)$, we can make $y_1$ and $y_2$ arbitrarily close to zero, and, by equations (\ref{eq:x_1}), (\ref{eq:z_1}), (\ref{eq:x_2}) and (\ref{eq:z_2}), make $x_1,x_2$ close to $1+\zeta$ and  $z_1,z_2$ close to $z^-$. This means that $Q_1$ and $Q_2$ are close to $M^-(0,1+\zeta,z^-)$. \qed
\subsubsection{Quasi-transverse intersection $W^s(Q)\cap W^u(P)$}\label{quasitrans1}
\par{}
The next step is to find integers $j_1$ and $j_2$, and parameter values of $\zeta$ and $\rho$ such that the index-2 periodic point $Q$ given by Lemma \ref{lem:p2i2} satisfies $W^s(Q)\cap W^u(P)\neq \emptyset$. This intersection is quasi-transverse. Indeed, we are going to consider the intersection of two local manifolds $W^s_{loc}(Q)$ and $W^u_{loc}(P)$. Note that $W^s_{loc}(Q)$ is a leaf of the foliation $\mathcal{F}_1$ tangent to strong-stable directions (i.e. $z$-directions), and $W^u_{loc}(P)$ is tangent to the center unstable direction (i.e. $(y,x)$-directions). Therefore, for the intersection point $x$, we have $\mathcal{T}_x W^s_{loc}(Q) \cap \mathcal{T}_x W^u_{loc}(P)=\{0\}$, which gives the quasi-transversality.
\par{}
Note that, by changing $\zeta$, one can move $M^-$. Consequently, we can control the position of the point $Q$ since it can be chosen arbitrarily close to $M^-$.  We will also change $\rho$ at the same time to ensure that $Q$ remains a index-2 periodic point.
\par{}
The following result holds.
\begin{lem}\label{lem:nontransverse}
For any $P\in\{P^+_k\}$ and any given $\rho^*\in(0,\dfrac{1}{2})$, there exists a sequence $\{\zeta_{i},\rho_i\}$ where $\zeta_i \to 0$ and $\rho_i\to\rho^*$ as $i\to +\infty$ such that, for each pair $(\zeta_{i},\rho_i)$ of parameter values, the corresponding system $X_{\zeta_{i},\rho_i}$ has a periodic point $Q_i$ of index 2 and period 2, where $Q_i \to M^-$ as $i\to +\infty$, and its stable manifold $W^s(Q_i)$ intersects the unstable manifold $W^u(P)$. 
\end{lem}
\par{}
\noindent \textit{Proof.} We consider an index-2 periodic orbit $\{Q=Q^1(y_1,x_1,z_1),Q^2(y_2,x_2,z_2)\}$ given by Lemma \ref{lem:p2i2}. In order to find the desired intersection, we need formulas for $W^u(P)$ and $W^s(Q)$. We pick an arbitrary point $P\in\{P_k\}$. Recall that those points $P_k$ are found in the proof of Lemma \ref{lem:WsP_1}. Let $P$ has the coordinates $(y_p,x_p,z_p)$. By taking a vertical line joining $P$ and a point on $\{y=0\}$ and iterating it, one can check that the local unstable manifold $W^u_{loc}(P)$ is spiral-like and winds onto $M^+$, which is given by
\begin{equation}\label{eq:WuP}
\begin{array}{rcl}
y&=&Ax_pt^{\rho}\cos(\omega\ln\dfrac{1}{t}+\eta)+o(t^\rho) \,,\\[7pt]
x&=&1+A_1x_pt^{\rho}\cos(\omega\ln\dfrac{1}{t}+\eta_1)+o(t^\rho) \,, \\[7pt]
z&=&z^+ +
                \begin{pmatrix}
                 A_2x_p t^\rho\cos{(\omega \ln {\dfrac{1}{t}}+\eta_2)}+o(t^\rho) \\
                 \cdots \\
                  A_{n-2}x_p t^\rho\cos(\omega \ln {\dfrac{1}{t}}+\eta_{n-2})+o(t^\rho) 
                \end{pmatrix} \,,
\end{array}
\end{equation}
\noindent where $t\in(0,y_p)$.
\par{}
In order to find the local stable manifold $W^s_{loc}(Q)$, we remind that there exists a absolutely continuous foliation $\mathcal{F}_1$ on $\Pi$ ,and $W^s_{loc}(Q)$ and $W^{ss}_{loc}(M^-)$ are leaves of $\mathcal{F}_1$ (see discussion after the non-degeneracy condition). By choosing $j_1$ and $j_2$ sufficiently large in (\ref{eq:p2i2}), we can make $Q$ arbitrarily close to $M^-$, and, therefore, $W^s_{loc}(Q)$ is arbitrarily close to $W^{ss}_{loc}(M^-)$ by the continuity of the foliation $\mathcal{F}_1$. Let us first find the formula for $W^{ss}_{loc}(M^-)$, and then we can obtain the formula for $W^{s}_{loc}(Q)$ by adding some small corrections. In this part, we write coordinates $x$ with its subscript as introduced in the beginning of this paper, i.e. $x=(x_1,x_2)$. Recall that the leaves of the foliation $\mathcal{F}_1$ on the cross-section $\Pi$ are obtained as the intersections of $\Pi$ with the orbits of the leaves of the foliation $\mathcal{F}_0$ by the flow. Since
the leaves of $\mathcal{F}_0$ on $W^s_{loc}(O)$ take the form $\{y=0,x_1=c_1,x_2=c_2\}$ where $c_1$ and $c_2$ are constants, and the cross-section $\Pi$ is a small piece of $\{x_2=0\}$, the leaves of $\mathcal{F}_1$ on $\Pi\cap W^s_{loc}(O)$ take the form $\{y=0,x_1=const\}$. Note that $W^{ss}_{loc}(M^-)\subset\Pi\cap W^s(O)$ and the $x_1$-coordinate of $M^-$ is $1+\zeta$. This implies that the local strong-stable manifold $W^{ss}_{loc}(M^-)$ is given by $\{y=0,x_1=1+\zeta\}$. Thus, the local stable manifold $W^s_{loc}(Q)$ has the form (we drop the subscript of $x_1$ again from now on)
\begin{equation}\label{eq:pp1_wsq}
\begin{array}{rcl}
y&=&f_1(z,\zeta,\rho,j_1,j_2)  \,,\\[7pt]
x&=&1+\zeta+f_2(z,\zeta,\rho,j_1,j_2) \,,
\end{array}
\end{equation}
\noindent where $f_1,f_2\to 0$ as $j_1,j_2\to +\infty$. 
\par{}
The intersection points of $W^{s}(Q)$ with $W^u(P)$ are given by equations (\ref{eq:WuP}) and (\ref{eq:pp1_wsq}). By noting $x_p=1+O(y_p^\rho)$ (since $P$ is a fixed point of $T_2$) and $z=z^+ + O(t^\rho)$ in (\ref{eq:WuP}), finding the intersection $W^{s}(Q)\cap W^u(P)$ is equivalent to solving the equations
\begin{eqnarray}
0&=&At^\rho\cos(\omega\ln\dfrac{1}{t}+\eta)+o(t^\rho)-f_1(z^+ + O(t^\rho),\zeta,\rho,j_1,j_2)\,,\label{eq:pp1_3a}\\[7pt]
\zeta&=&A_1t^{\rho}\cos(\omega\ln\dfrac{1}{t}+\eta_1)+o(t^\rho)-f_2(z^+ + O(t^\rho),\zeta,\rho,j_1,j_2) \,. \label{eq:pp1_3b}
\end{eqnarray}
\noindent The equation (\ref{eq:pp1_3a}) gives a countable set $\{t_i\}$ of $t$ values where $t_i\to 0$ as $i\to +\infty$, and $t_i$ are functions of $j_1,j_2,\rho$ and $\zeta$. We plug $t_i$ into (\ref{eq:pp1_3b}) and get
\begin{equation}\label{eq:pp1_4}
\zeta - f_3(\rho,\zeta,i,j_1,j_2)= 0 \,,
\end{equation}
\noindent where $f_3$ is continuous and tend to 0 as $i,j_1,j_2\to +\infty$. Next, we pick an arbitrary sequence $\{(j_1^i,j_2^i)\}$ such that $j_1^i,j_2^i \to +\infty$ and $\dfrac{j_1^i}{j_2^i}\to \rho^*$  as $i\to +\infty$. Note that, for such sequence, we can obtain a formula for $\rho$ from equation (\ref{eq:p2i2}) by implicit function theorem. Indeed, by plugging $(j_1^i,j_2^i)$ into (\ref{eq:p2i2}) and sorting the terms, we have
\begin{equation*}
\rho=\dfrac{j_2^i}{j_1^i}+\dfrac{\Psi(\zeta,\rho,j_1^i,j_2^i,c)}{j_1^i}.
\end{equation*}
\noindent Since the second term in the RHS of the above equation tends to zero as $i$ tends to positive infinity, we can, by implicit function theorem, rewrite this equation as
\begin{equation}\label{eq:pp1_4.1}
\rho=\dfrac{j_2^i}{j_1^i}+\dfrac{\hat{\Psi}(\zeta,j_1^i,j_2^i,c)}{j_1^i},
\end{equation}
\noindent where $\hat{\Psi}$ is continuous in $\zeta$ and $c$, and is uniformly bounded. We now let $g(\zeta,j_1^i,j_2^i):=\dfrac{\Psi(\zeta,j_1^i,j_2^i,c^*)}{j_1^i}$, where $c^*$ can be any value in $(-1,1)$. Then, by Lemma \ref{lem:index_2} and equation (\ref{eq:pp1_4}), the index-2 point $Q$ whose stable manifold $W^s(Q)$ intersects the unstable manifold $W^u(P)$ of the index-1 fixed point $P$ can be found by solving the following system of equations
\begin{equation}\label{eq:pp1_5} 
\begin{array}{rcl}
\rho- \dfrac{j_2^i}{j_1^i} - g(\zeta,j_1^i,j_2^i)& =& 0,  \\[10pt]
\zeta - f_3(\rho,\zeta,i,j_1^i,j_2^i)&= &0    .
\end{array}  
\end{equation}
\noindent By plugging the first equation of (\ref{eq:pp1_5}) into the second one, we have 
\begin{equation}\label{eq:pp1_6} 
\zeta = f_3\Big(\dfrac{j_2^i}{j_1^i} + g(\zeta,j_1^i,j_2^i),\zeta,i,j^i_1,j^i_2\Big) .
\end{equation}
\noindent Since the RHS of equation (\ref{eq:pp1_6}) is continuous and tends to zero as $i$ tends to positive infinity, for each sufficiently large $i$, we can find parameter values $\zeta_i$ satisfying (\ref{eq:pp1_6}), and then, from \eqref{eq:pp1_5}, find $\rho_i$ satisfying (\ref{eq:p2i2}) where $\zeta_i \to 0$ and $\rho_i \to \rho^*$ as $i \to +\infty$. This means that, in the corresponding system $X_{\rho_i,\zeta_i}$, the point $Q$ has period 2 and index 2 and its stable manifold $W^s(Q)$ intersects the unstable manifold $W^u(P)$ of the index-1 fixed point $P$. \qed
\subsubsection{Transverse intersection $W^u(Q)\cap W^s(P)$}\label{prf:transverse1}
Lemma \ref{lem:nontransverse} implies that a heterodimensional cycle will be created if, for any pair $(\zeta_{i},\rho_i)$ in this lemma, the corresponding index-2 periodic point $Q$ also satisfy $W^u(Q) \cap W^s(P)\neq\emptyset$. We now prove that this intersection exists and it is transverse.  We note the following result on the unstable manifold.
\begin{lem}\label{lem:WuQ}
The unstable manifold $W^u(Q)$ of the orbit of an index-2 periodic point $Q$ intersects $\Pi_0=\{y=0\}$ transversely.
\end{lem}
\noindent \textit{Proof.} Consider the map $\tilde{T}\equiv(\tilde{T}_1,\tilde{T}_2)$ obtained from the Poincaré map $T\equiv(T_1,T_2)$ by taking quotient along leaves of the strong-stable foliation on $\Pi$. This map $\tilde{T}$ acts on the 2-dimensional surface $\tilde{\Pi}=\Pi\cap \{z=z^*\}$ where $\|z^*\|<\delta$. We call this surface $\tilde{\Pi}$ a quotient cross-section. More specifically, for a region $V\subset\tilde{\Pi}$, its image $\tilde{T}(V)$  is the projection of $T(V)$ onto $\tilde{\Pi}$ along the leaves of the strong-stable foliation $\mathcal{F}_1$.
\par{}
For any region $V\subset\tilde{\Pi}$ such that $T(V)\subset\Pi$, we have that
\begin{equation}\label{qm5}
qS(V)<S(\tilde{T}(V)) \,,
\end{equation}
\noindent where $S$ denotes the area and $q>1$. We now prove this inequality. Let $\tilde{\Pi}_1=\tilde{\Pi}\cap\{y\geqslant0\}$ and $\tilde{\Pi}_2=\tilde{\Pi}\cap\{y\leqslant0\}$. We assume that $S(V\cap \tilde{\Pi}_1)>S(V\cap \tilde{\Pi}_2)$ and proceed by considering the region $V\cap \tilde{\Pi}_1$. We have $S(V\cap \tilde{\Pi}_1)>\dfrac{1}{2}S(V)$. Let us first look at the equations on coordinates $y$ and $x$ in the formula (\ref{eq:setting_7}) of the map $T_1$, i.e. 
\begin{equation}\label{qm1}
\left\{\begin{array}{l}
                \bar{y}=Axy^\rho\cos{(\omega \ln {\dfrac{1}{y}}+\eta)}+o(y^\rho) \\[7pt]
                \bar{x}=1+A_1xy^\rho\cos{(\omega \ln {\dfrac{1}{y}}+\eta_1)}+o(y^\rho)
                \end{array}\right..\
\end{equation}
\noindent We have that 
\begin{equation}\label{qm3}
\left|
\dfrac{\partial(\bar{y},\bar{x})}{\partial y\partial x}\right|=-\omega AA_1\sin(\eta_1-\eta)|y|^{2\rho-1}+o(|y|^{2\rho-1})\,,
\end{equation}
\noindent where $AA_1\sin(\eta_1-\eta)$ is non-zero by the non-degeneracy condition (see (\ref{eq:nondege})). The determinant (\ref{qm3}) is much greater than one since $y$ is small and $\rho<\dfrac{1}{2}$. This implies that the projection of the 2-dimensional area of $T_1(V\cap\tilde{\Pi}_1)$ on the $(x,y)$-plane is much larger than the area of $V\cap\tilde{\Pi}_1$ and, therefore, the area of $V$. it follows that the image $T(V)$ (which contains $T_1(V\cap\tilde{\Pi}_1)$) has an area much larger than that of $V$. Note that the derivatives $\dfrac{\partial \bar{z}}{\partial y}$ and $\dfrac{\partial \bar{z}}{\partial y}$ from the formula (\ref{eq:setting_7}) are so small that the angle between the image $T(V)$ and the horizontal surface is also small. Therefore, the leaves of the foliation $\mathcal{F}_1$ are transverse to $T(V)$. It follows that projecting along those leaves changes the area of $T(V)$ by a factor that is finite and bounded away from zero. Therefore, the image $\tilde{T}(V)$ has an area much larger than $V$. The inequality (\ref{qm5}) follows.
\par{}
Let $\tilde{Q}$ be the projection of $Q$ on $\tilde{\Pi}$ along the leaves of the strong-stable foliation $\mathcal{F}_1$. Then the point $\tilde{Q}$ is a completely unstable periodic point of $\tilde{T}$. Let $U$ be a small neighborhood containing $\tilde{Q}$. We now show that there exists some $i$ such that the image $\tilde{T}^{(i)}(U)$ intersects $\tilde{\Pi}_0:=\Pi_0\cap\{z=z^*\}$ transversely. 
\par{}
We start by claiming that there are infinitely many pre-images of $\tilde{\Pi}_0$ under $\tilde{T}$ on $\tilde{\Pi}$ and they are nearly horizontal lines crossing $\tilde{\Pi}$. Let us first consider the pre-images of $\Pi_0$ under $T$, which are surfaces in $\Pi$. By the formulas (\ref{eq:setting_7}) and (\ref{eq:setting_8}) for the map $T$, these surfaces $T_1^{(-1)}(\Pi_0)$ and $T_2^{(-1)}(\Pi_0)$ are given by
\begin{equation}
0=Axy^\rho\cos{(\omega \ln {\dfrac{1}{y}}+\eta)}+o(y^\rho)\,,
\end{equation}
\noindent and 
\begin{equation}
0=B(1+\zeta)xy^\rho\cos{(\omega \ln {\dfrac{1}{|y|}}+\theta)}+o(|y|^\rho)\,,
\end{equation}
\noindent which, by the transformation (\ref{eq:pp1_3}), give
\begin{equation}\label{eq:yk1}
y_k^1=\exp\bigg(\dfrac{-(2k+1)\pi+2\eta}{2\omega}\bigg)+o(1)_{k\to+\infty}\,,
\end{equation}
\noindent and
\begin{equation}\label{eq:yk2}
y_k^2=-\exp\bigg(\dfrac{-(2k+1)\pi+2\theta}{2\omega}\bigg)+o(1)_{k\to+\infty}\,.
\end{equation}
These surfaces $y^i_k$ $(i=1,2)$ with any sufficiently large $k$ are the pre-images of $\{y=0\}$ under $T_i$. Those pre-images are pieces of $W^s(O)\cap\Pi$ which consists of leaves of the foliation $\mathcal{F}_1$. Note that $\tilde{\Pi}$ is transverse to those leaves. When we project $y_k^i$ onto $\tilde{\Pi}$ along the leaves, we get curves $l^i_k$, which are pre-images of $\{y=0\}$ under $\tilde{T}_i$. The claim follows.
\par{}
Note that we can choose $Q$ sufficiently close to $\Pi_0$, such that the orbit of $\tilde{Q}$ as well as the neighborhood $U$ are inside a region bounded by $\{x=1+\delta\}$, $\{x=1-\delta\}$, $l^2_{k_1}$ and $l^2_{k_2}$, for some $k_1$ and $k_2$. Now let us iterate the neighborhood $U$. On one hand, by (\ref{qm5}), the area of $\tilde{T}^{(i)}(U)$ is expanding as the number $i$ increases. On the other hand, from (\ref{eq:setting_7}) and (\ref{eq:setting_8}), the coordinate $x$ is uniformly bounded when $y$ is sufficiently small, and, therefore, the iterate $\tilde{T}^{(i)}(U)$ cannot intersect the boundaries $\{x=1+\delta\}$ and $\{x=1-\delta\}$. It follows that there exists some integer $i$ such that the image $\tilde{T}^{(i)}(U)$ intersects transversely either $\tilde{\Pi}_0$ or one of the boundaries $l^2_{k_1}$ and $l^2_{k_2}$. For the latter case, the next iterate $\tilde{T}^{(i+1)}(U)$ intersects $\tilde{\Pi}_0$ transversely. 
\par{}
Let us now consider a disc $U_0$ which is centered at $Q$ and paralleled to $(x,y)$-plane. Let $U$ be the projection of $U_0$ on $\tilde{\Pi}$ along the leaves. By the result above, we have that $\tilde{T}^{(i)}(U)$ intersects $\{y=0\}$ transversely for some $i$. By the way how we define the map $\tilde{T}$, we have that $T^{(i)}(U_0)$ intersects $\{y=0\}$ transversely. Since the unstable manifold $W^u(Q)$ is obtained by taking limit of the iterates $T^{(n)}(U_0)$, we have that $W^u(Q)$ intersects $\{y=0\}$ transversely. \qed
\par
\begin{figure}[!h]
\centering
\includegraphics[width=0.55\textwidth]{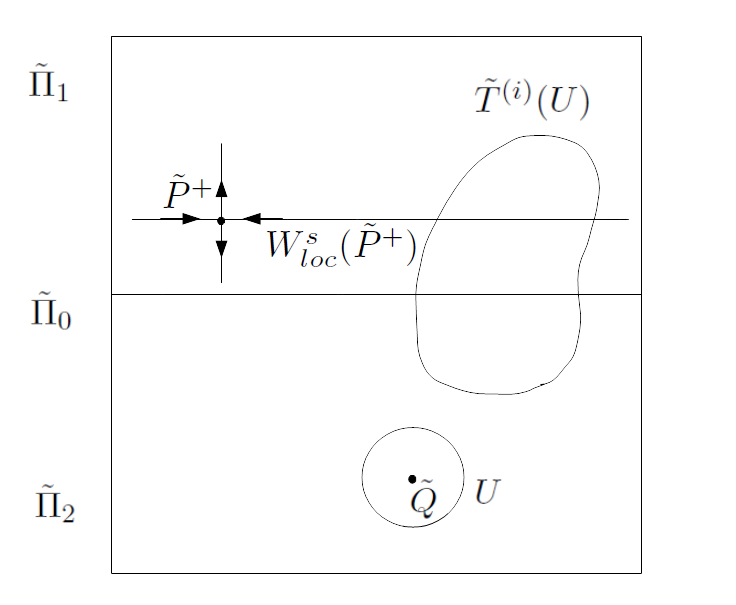}
\caption{A schematic picture for $W^u(\tilde{Q})\cap W^s(\tilde{P}^+)$, where the point $\tilde{P}^+$ corresponds to the point $P^+\in\{P^+_k\}$.}
\label{fig:transintersection1}
\end{figure}

{}
By combining Lemma \ref{lem:WuQ} and Lemma \ref{lem:WsP_1}, we have that there exists a point $P^+\in\{P_k^+\}$ such that $W^u(Q)$ intersects $W^s(P^+)$ transversely (see figure \ref{fig:transintersection1}). If the $y$-coordinate of $P$ is smaller than that of $P^+$, then, by Lemma \ref{lem:WsP_1}, $W^s(P)$ will be below $W^s(P^+)$, which implies $W^u(Q)\cap W^s(P)\neq\emptyset$. We now consider the case where the $y$-coordinate of $P$ is larger than that of $P^+$ (i.e. in the sequence $\{P^+_k\}$, $P^+$ has a subscript larger than that of $P$). If we can show $W^u(P^+)\cap W^s(P)\neq\emptyset$, then we will obtain $W^u(Q)\cap W^s(P)\neq\emptyset$ by $\lambda$-lemma.
\par
\begin{figure}[!h]
\centering
\includegraphics[width=0.95\textwidth]{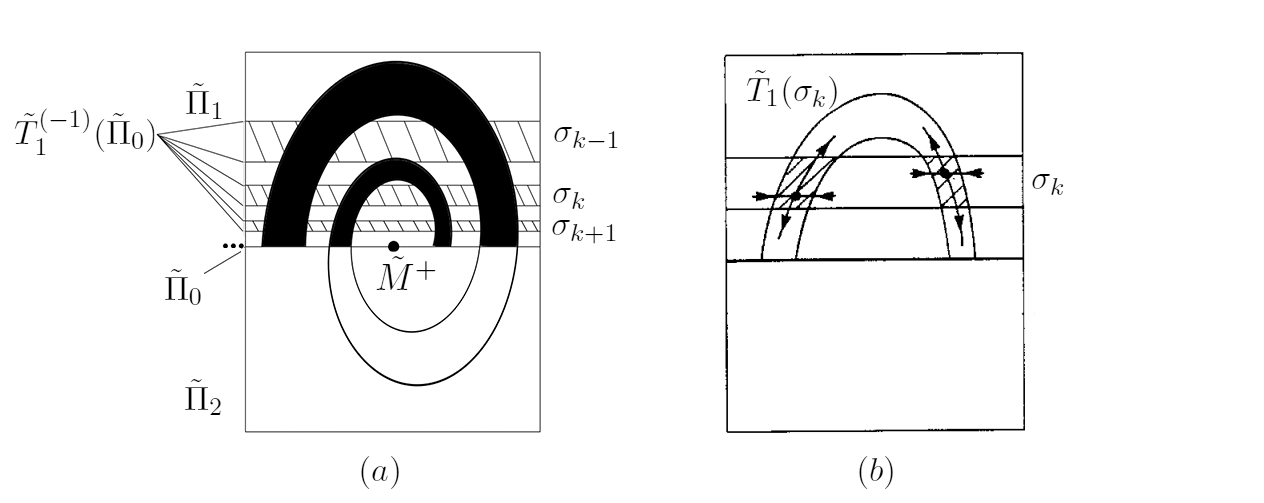}
\caption{Here we show the horseshoes in the quotient cross-section $\tilde{\Pi}$. We still denote by $\sigma_k$ the intersection $\sigma_k\cap\tilde{\Pi}$ and it is bounded by curves $l^i_{2k-1}$ and $l^i_{2k}$ defined under the equation (\ref{eq:yk1}). The two black-filled curls (from inside to outside) are images of $\sigma_{k+1}$ and $\sigma_k$, respectively.}
\label{fig:transintersection2}
\end{figure}

\begin{figure}[!h]
\centering
\includegraphics[width=0.60\textwidth]{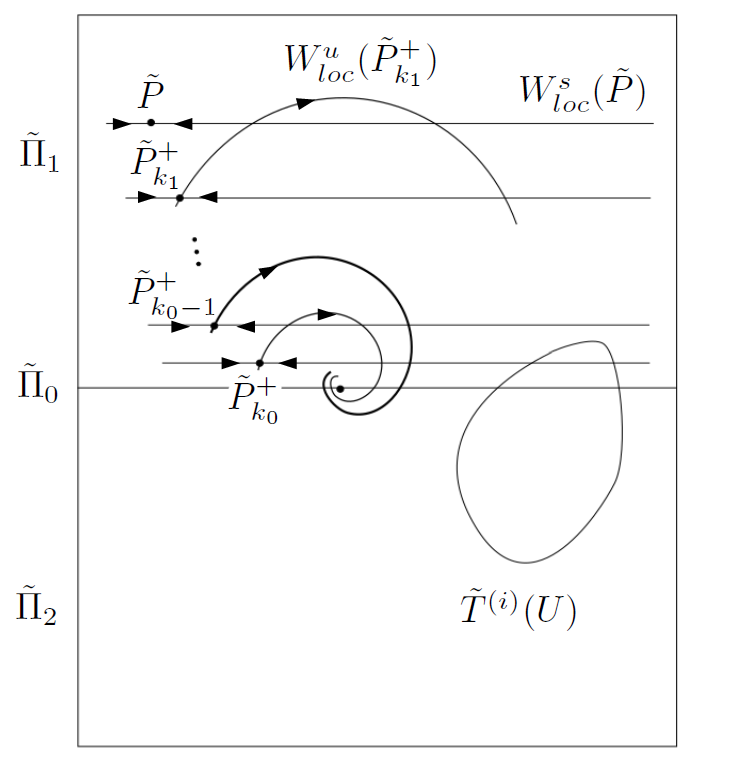}
\caption{In the quotient cross-section $\tilde{\Pi}$, the local stable manifolds $W^s_{loc}(\tilde{P}_k^+)$ are horizontal lines and the local unstable manifolds $W^u_{loc}(\tilde{P}_k^+)$ are spirals winding onto $\tilde{M}^+$. The image $\tilde{T}^{(i)}(U)$ contains a piece of $W^u(\tilde{Q})$ that intersects $W^s(P^+_{k_0})$. The points $P_k^+$ are homoclinically related, and, especially, we have $W^u(P^+_{k_0})\cap W^s(P)\neq\emptyset$, where the point $P^+_{k_1}$ is the one just after $P$ in the sequence $\{P^+_k\}$. }
\label{fig:transintersection3}
\end{figure}

{}
Let us now establish the heteroclinic connection $W^u(P^+)\cap W^s(P)\neq\emptyset$. Denote by $\sigma_k$ the regions in $\Pi_1$ bounded by surfaces $\{y=y^1_{2k-1}\}$ and $\{y=y^1_{2k}\}$ given by formulas (\ref{eq:yk1}) and (\ref{eq:yk2}). One can easily check that the images $T_1(\sigma_k)$ belong to $\Pi_1$ and points outside $\sigma_k$ are mapped in to $\Pi_2$. By Shilnikov theorem there are, at $\rho<1$, infinitely many horseshoes in the cross-section $\Pi_1$, each of which corresponds to a region $\sigma_k$ and its image $T_1(\sigma_k)$, and there are two fixed points in each region (see figure \ref{fig:transintersection2}). Obviously, $W^u(P_i^+)\cap W^s(P_j^+)\neq\emptyset$ if $T_1(\sigma_i)$ intersects $\sigma_j\neq\emptyset$ properly. Here "properly" means that the intersection $T_1(\sigma_i)\cap\sigma_j$ is connected and the map $T_1|_{\sigma_i\cap T^{-1}_1(\sigma_j)}$ is a saddle map in the sense of \citep{sh67}. It is shown in \citep{sh70} that we have $T_1(\sigma_i)\cap \sigma_j\neq\emptyset$ if $j>\rho'i$ where $\rho'>\rho$ can be chosen arbitrarily close to $\rho$. We now pick $P^+=P_{k_0}^+\in\{P_k^+\}$ such that $W^u(Q)\cap W^s(P_{k_0}^+)\neq\emptyset$. We have that $W^s(P_{k_0-1}^+)\cap W^u(P_{k_0}^+)\neq\emptyset$ as long as $\rho'k_0<(k_0+1)$. Similarly, we obtain $W^s(P_{k_0-2}^+)\cap W^u(P_{k_0-1}^+)\neq\emptyset$ if $\rho'(k_0-1)<(k_0-2)$. Indeed, we can repeat this procedure until we arrive at $k=1$ since we can choose $\rho<\rho'<\dfrac{1}{2}$. Therefore, we find the heteroclinic intersections $W^u(P^+_{k_0})\cap W^s(P^+_k)$ where $1\leqslant k<k_0$. By noting that we assumed that the point $P$ is from $\{P_k^+\}$ where $1\leqslant k<k_0$, we obtain $W^u(P^+_{k_0})\cap W^s(P)\neq\emptyset$ (see figure \ref{fig:transintersection3}). Therefore, the existence of the transverse intersection $W^u(Q)\cap W^s(P)$ is proved.
\par{}
This transverse intersection $W^u(Q)\cap W^s(P)$ along with the quasi-transverse intersection $W^s(Q)\cap W^u(P)$ given by Lemma \ref{lem:nontransverse} give rise to a heterodimensional cycle of the map $T$ corresponding to the system $X_{\zeta_i,\rho_i}$. The theorem is proved. 
\subsection{Proof of Theorem \ref{thm2}}
We will find a heterodimensional cycle of the Poincaré map $T$ related to an index-1 fixed point $P$ and an index-2 period-3 point $Q$ by the procedure similar to that used in the proof of Theorem \ref{thm1}. The major difference is that, instead of directly finding an index-2 periodic point such that $W^s(Q)\cap W^u(P)\neq\emptyset$, we now find a point $Q$ of period 3 and index 2 with the property that the point $M^+$ falls onto the stable manifold $W^s(Q)$.
\subsubsection{Construction of the Poincaré map $T$}
Let us use the same cross-section $\Pi$ defined in the proof of Theorem \ref{thm1}. Recall that the intersection points $M^+$ and $M^-$ of the two loops $\Gamma^+$ and $\Gamma^-$ with $\Pi$ have coordinates $(y^+,x^+,z^+)$ and $(y^-,x^-,z^-)$. In Theorem \ref{thm2}, we will spilt the two homoclinic loops by using two parameter $\mu_1=y^+$ and $\mu_2=y^-$. In particular, we will consider the case where $\mu_1=-\mu_2$. The Poincaré map $T$ is now slightly different from that in proof of Theorem \ref{thm1}. The local maps $T_{loc_1}$ and $T_{loc_2}$ are the same. For the global maps $T_{glo_1}$ and $T_{glo_1}$, we add $\mu$ and $-\mu$ to the right hand side of the second equations in (\ref{eq:setting_4}) and (\ref{eq:setting_5}), respectively. After additionally replacing $\dfrac{\mu}{d}$ by $\mu$ in the compositions $T_{glo_1}\circ T_{loc_1}$ and $T_{glo_2}\circ T_{loc_2}$, we obtain the Poincaré map $T\equiv(T_1,T_2)$ given by
\begin{equation}\label{eq:setting_7a}
T_1:\quad\left\{\begin{array}{l}

                \bar{y}=\mu+Axy^\rho\cos{(\omega \ln {\dfrac{1}{y}}+\eta)}+o(y^\rho) \\[7pt]
                \bar{x}=1+A_1xy^\rho\cos{(\omega \ln {\dfrac{1}{y}}+\eta_1)}+o(y^\rho) \\[7pt]
                \bar{z}=z^+ +
                \begin{pmatrix}
                 A_2xy^\rho\cos{(\omega \ln {\dfrac{1}{y}}+\eta_2)}+o(y^\rho) \\
                 \cdots \\
                  A_{n-2}xy^\rho\cos(\omega \ln {\dfrac{1}{y}}+\eta_{n-2})+o(y^\rho) 
                \end{pmatrix}                
                \end{array}\right.\,,
\end{equation}
\noindent and
\begin{equation}\label{eq:setting_8a}
T_2:\quad\left\{\begin{array}{l}
                \bar{y}=-\mu-Bx|y|^\rho\cos{(\omega \ln {\dfrac{1}{|y|}}+\theta)}+o(|y|^\rho) \\[7pt]
                \bar{x}=1+\zeta+B_1x|y|^\rho\cos{(\omega \ln {\dfrac{1}{|y|}}+\theta_1)}+o(|y|^\rho) \\[7pt]
                \bar{z}=z^- + 
                \begin{pmatrix}
                 B_2x|y|^\rho\cos{(\omega \ln {\dfrac{1}{|y|}}+\theta_2)}+o(|y|^\rho) \\
                 \cdots \\
                  B_{n-2}x|y|^\rho\cos(\omega \ln {\dfrac{1}{|y|}}+\theta_{n-2})+o(|y|^\rho) 
                \end{pmatrix}     
                 \end{array}\right.\,,
\end{equation}\\
\noindent where $z \in \mathbb{R}^{n-3}$ and the coefficients are defined in the same way as those in the Poincaré map in Theorem \ref{thm1} given by (\ref{eq:setting_7}) and (\ref{eq:setting_8}).
\par{}
We now consider a 5-parameter family $X_{\mu,\zeta,\rho,u,v}$ of perturbed systems so $x^\pm,y^\pm$ and $z^\pm$ are smooth functions of parameters. The parameters $u$ and $v$ are smooth functions of coefficients of the system given by
\begin{equation*}
u=\frac{1}{2\pi}
\bigg(\omega\log{(B\sin{\varphi})}-\theta+\rho\theta-\rho\bigg(\frac{3\pi}{2}+\varphi\bigg)+\frac{\pi}{2}\bigg) \,,
\end{equation*}
\noindent and
\begin{equation*}
v=\frac{1}{2\pi}
\bigg(\omega\log{\bigg(\frac{1}{2}B\sin{(\theta_1-\theta)}\bigg)}-\theta+\rho\theta-\rho\bigg(\frac{\pi}{2}+\theta-\theta_1\bigg)+\frac{3\pi}{2}+\varphi\bigg) \,,
\end{equation*}
\noindent where $\varphi=\arctan\dfrac{\rho}{\omega}$. The reason for defining those two parameters can be seen in the proof of Lemma \ref{lem:nontransverse2} in Section \ref{prf:i2QM}. 
\subsubsection{Existence of the index-1 fixed point $P$}
Each point of $\{P^+_k\}$ and $\{P^-_k\}$ mentioned in the proof of Theorem \ref{thm1} remains an index-1 saddle fixed point under sufficiently small perturbations (the closer the point is to $W^s(O)$, the smaller the perturbation must be). We now pick a point $P$ from the set $\{P^+_k\}\cup\{P^-_k\}$ and consider perturbations under which $P$ is still an index-1 saddle fixed point (i.e. we choose $\mu$ sufficiently small). In what follows, we consider $P\in\{P^+_k\}$. We remark here that if we use a point $P\in\{P^-_k\}$, we can sill find a heterodimensional cycle in a similar way, and the difference is that the functions defining $u$ and $v$ will change since the coefficients in $T_1$, instead of $T_2$, are involved. 
\par{}
We now state a lemma on points in $\{P^+_k\}$ under small perturbation.  We will show later that the point $P$ given by this lemma can be the desired index-1 point to create a heterodimensional cycle.
\begin{lem}\label{lem:WsP_2}
For any $\mu$ sufficiently close to 0, there exists a point $P\in\{P_k^+\}$ remaining a saddle fixed point of $T_1$ in the system $X_{\mu,\zeta,\rho,u,v}$ such that the stable manifold $W^s(P)$ is the graph of a function of coordinates $x$ and $z$ defined for all $x$ and $z$ values in $\Pi$ and is bounded by $\{y=|\mu|\}$ and $\{y=0\}$.
\end{lem}
\par{} 
\noindent \textit{Proof.} The proof of this lemma is similar to that of Lemma \ref{lem:WsP_1}. At $\mu=0$, the fixed points $P^+_k$ are given by
\begin{equation}\label{eq:WsP2_1}
\begin{array}{rcl}
y_k&=&C\exp\Big(\dfrac{-\pi k}{\omega}\Big)+o\Big(\exp\Big(\dfrac{-\pi k}{\omega}\Big)\Big)\, ,\\[10pt]
x_k&=&1+o\Big(\exp\Big(\dfrac{-\pi k}{\omega}\Big)\Big)\, ,\\[10pt]
z_k&=&z^+ +o\Big(\exp\Big(\dfrac{-\pi k}{\omega}\Big)\Big)\, ,
\end{array}
\end{equation}
\noindent where $y_k,x_k,z_k$ are the coordinates of $P^+_k$, $C=\exp\Big(\dfrac{2\eta-\pi}{2\omega}\Big)$, and $k$ is any positive integer greater than some sufficiently large $K$. Equations in (\ref{eq:WsP2_1}) hold at $\mu\neq 0$ if $\mu\exp(\dfrac{\pi\rho k}{\omega})$ is sufficiently close to 0. By the argument under (\ref{eq:WsP_2}) in the proof of Lemma \ref{lem:WsP_1}, we have that $P^+_k$ is also of index 1 for sufficiently small values of $\mu$. 
\par{}
We now consider the inverse image under $T_1$ of a small piece of the surface $\{y=y_k\}$ containing $P^+_k$. By (\ref{eq:setting_7a}), we have
\begin{equation}\label{eq:WsP2_2}
\sin\Big(\dfrac{\pi}{2}-\theta-\omega\ln\dfrac{1}{y}\Big)=\dfrac{1}{xA}\Big(\dfrac{y_k}{y^\rho}-\dfrac{\mu}{y^\rho}+o(1)_{y\to 0}\Big)\, ,
\end{equation}
\noindent where $(y,x)$ are coordinates of the points in the inverse image ($z$ coordinates are in the small term) and $x$ is bounded since the small cross-section $\Pi$ is bounded. We have following equation if $y,\dfrac{y_k}{y^\rho}$ and $\dfrac{|\mu|}{y^\rho}$ are sufficiently small:
\begin{equation}\label{eq:WsP2_3}
\dfrac{\pi}{2}-\theta-\omega\ln\dfrac{1}{y}=\dfrac{1}{xA}\Big(\dfrac{y_k}{y^\rho}-\dfrac{\mu}{y^\rho}+o(1)_{y\to 0}\Big)+m\pi\quad ,m=0,\pm1,\pm2,\ldots \quad,
\end{equation}
\noindent which, by noting that the surface contains $P^+_k$, leads to
\begin{equation}\label{eq:WsP2_4}
y=C\exp(\dfrac{-\pi k}{\omega})+o(1)_{k\to +\infty}\, .
\end{equation}
\noindent Formula (\ref{eq:WsP2_4}) has the same form as (\ref{eq:WsP_5}), and it is valid for all values of $x,z$, where $(y,x,z)\in\Pi$, if $|\mu|\exp(\dfrac{\pi\rho k}{\omega}),y,\dfrac{y_k}{y^\rho}$ and $\dfrac{\mu}{y^\rho}$ are sufficiently small. From equations in (\ref{eq:WsP2_1}) and formula (\ref{eq:WsP2_4}), this requirement is equivalent to that $k$ is sufficiently large and $|\mu|\exp(\dfrac{\pi\rho k}{\omega})$ is sufficiently small. This can be satisfied since $\rho<1$ and we can choose sufficiently small $\mu$ and sufficiently large $k$ independently. Especially, $\mu$ and $k$ can be chosen such that $|\mu|\gg y_k$. Indeed, by letting $\mu=\exp(\dfrac{-2\pi j-\xi_\mu+\eta}{\omega})$, to obtain $|\mu|\gg y_k$ (while $|\mu|\exp(\dfrac{\pi\rho k}{\omega})$ is small) is equivalent to find $j$ and $k$ such that $\rho k \ll j \ll k$. One can check that the successive backward iterates of a small piece of the surface $y=y_k$ containing $P^+_k$ take the form as (\ref{eq:WsP2_3}), where the term $o(1)$ stays uniformly small. Since $W^s_{loc}(P)$ is the limit of a sequence of those iterates, $W^s_{loc}(P)$ is given by (\ref{eq:WsP2_3}). \qed
\subsubsection{An index-2 periodic point $Q$ with $M^+\in W^s(Q)$}\label{prf:i2QM}
Here we consider a periodic orbit of $T$ such that it not only has index 2 but also satisfies the property that the point $M^+$ falls onto its stable manifold. Such orbit allows for the emergence of a quasi-transverse intersection $W^s(Q)\cap W^u(P)$ after an arbitrarily small perturbation in $\zeta$. The following result holds.
\begin{spacing}{1.3}\begin{lem}\label{lem:nontransverse2}
Let $(\hat{\rho},\hat{u},\hat{v})$ be a triple such that $\hat{\rho}=\dfrac{p}{q}\in(\mathbb{Q}\cap(0,\dfrac{1}{2})),\hat{u}=\dfrac{p_1}{q}$ and $\hat{v}=\dfrac{p_2}{q}$, where $p,q$ are co-prime and $p_1,p_2$ are any integers. The triple $(\hat{\rho},\hat{u},\hat{v})$ corresponds to a sequence $\{(\mu_{j},\zeta_j,u_{j}, v_{j})\}$ accumulating on $(0,0,\hat{u},\hat{v})$ such that the map $T$ corresponding to the system $X_{\mu_j,\zeta_j,\hat{\rho},u_{j}, v_{j}}$ has a periodic point $Q$ of period 3 and index 2 satisfying that $M^+\in W^s(Q)$.
\end{lem}\end{spacing}
\par{}
We remark here that this lemma is for the case where the saddle index $\rho^*$ of the unperturbed system $X$ is rational; when it is not, we only need to do an arbitrarily small perturbation.
\begin{figure}[!h]
\begin{center}
\includegraphics[scale=0.6]{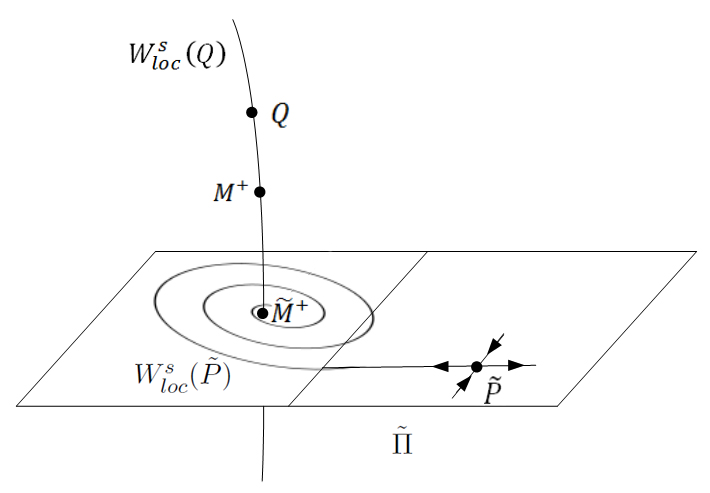}
\end{center}
\caption{The configuration of objects in $\tilde{\Pi}$ when $M^+\in W^{s}(Q)$.}
\label{fig:method2}
\end{figure}
\par{}
\noindent \textit{Proof of Lemma \ref{lem:nontransverse2}.} Consider now a periodic orbit $\{Q=Q_1(y_1,x_1,z_1),Q_2(y_2,x_2,z_2),Q_3(y_3,$ $x_3,$ $z_3)\}$ of period 3. We will show that there exist parameter values for which $Q_1$ has index 2, and the point $M^+$ lies on the local stable manifold $W^{s}_{loc}(Q_1)$ (figure \ref{fig:method2}). 
\begin{figure}[!h]
\begin{center}
\includegraphics[scale=0.75]{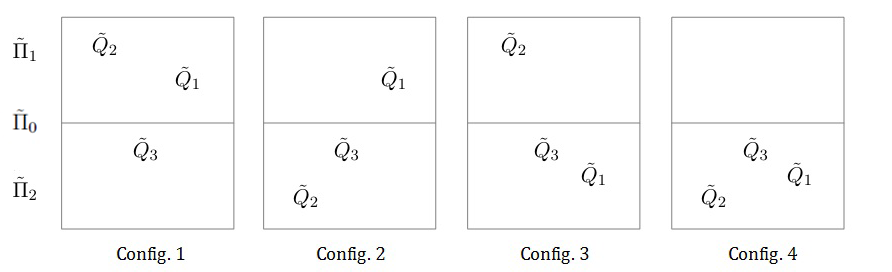}
\end{center}
\caption{The configurations of points $\tilde{Q_1},$ $\tilde{Q_2}$ and $\tilde{Q_3}$ on the quotient cross-section $\tilde{\Pi}$.}
\label{fig:config}
\end{figure}
\par{}
Depending on the sign of $\mu$ and the values of $\theta_1-\theta$, there are four logical possibilities of configurations of the points $Q_1,Q_2$ and $Q_3$ on $\Pi$, (see figure \ref{fig:config}), which are given by
\begin{equation}\label{eq:configurations}
\begin{array}{ll}
1.& Q_1=T_2(Q_3)\in\Pi_1\,,\quad Q_2=T_1(Q_1)\in\Pi_1\,,\quad Q_3=T_1(Q_2)\in\Pi_2\,; \\[5pt]
2.& Q_1=T_2(Q_3)\in\Pi_1\,,\quad Q_2=T_1(Q_1)\in\Pi_2\,,\quad Q_3=T_2(Q_2)\in\Pi_2\,, \\[5pt]
3.& Q_1=T_2(Q_3)\in\Pi_2\,,\quad Q_2=T_2(Q_1)\in\Pi_1\,,\quad Q_3=T_1(Q_2)\in\Pi_2\,; \\[5pt]
4.& Q_1=T_2(Q_3)\in\Pi_2\,,\quad Q_2=T_2(Q_1)\in\Pi_2\,,\quad Q_3=T_2(Q_2)\in\Pi_2\,. \\[5pt]
\end{array}
\end{equation}
\noindent Note that, for different configuration, the formulas for the parameters $u$ and $v$ will change but the same result in Lemma \ref{lem:nontransverse2} holds. 
\par{}
Here we only consider the case where $\mu<0$ and $(\theta_1-\theta)\in(0,\dfrac{\pi}{2}]$, i.e. the fourth configuration in (\ref{eq:configurations}).
\par{}
We first need a formula for the local stable manifold $W^{s}_{loc}(Q_1)$, which is a leaf of the strong-stable foliation $\mathcal{F}_1$. The leaves of $\mathcal{F}_1$ are given by the following lemma.
\begin{lem}\label{lem:WsQ}\begin{spacing}{1.2}
Let $M(y_0,x_0,z_0)$ be a point on $\Pi$ with $y_0$ sufficiently small. The local strong stable manifold $W^{ss}_{loc}(M)$ (i.e. the leaf of $\mathcal{F}_1$ through $M$) is the graph of the function
\begin{equation*}
(y,x)=(y_0+(z-z_0)a_1,\,\,x_0+(z-z_0)a_2)\,,
\end{equation*}
\noindent where $a_1=o(|y_0|)$ and $a_2=o(1)_{y_0\to0}$ are $(n-3)$-dimensional vectors whose components are certain functions of $y_0,x_0,z_0$ and the parameters $\varepsilon$ satisfying \end{spacing}
\begin{align*}
\dfrac{\partial^{i+k+l+n} o(|y_0|)}{\partial^i y_0 \,\partial^k x_0 \,\partial^l z_0\,\partial^n \varepsilon}&=o(|y_0|^{1-i})\quad i+k+l+n\leqslant(r-1)\,,\\[20pt]
\dfrac{\partial^{k+l+n} o(1)}{\partial^k x_0 \,\partial^l z_0\,\partial^n \varepsilon}&=o(1)_{y_0\to 0} \quad k+l+n\leqslant(r-1) \,. 
\end{align*}
\end{lem}
\noindent This result is lemma 4 in \citep{lt15} and we omit the proof here. The local stable manifold $W^{s}_{loc}(Q_1)$ is now given by
\begin{equation}\label{eq:WsQ}
\begin{array}{rcl}
y&=&y_1+(z-z_1)\delta_1  \,,\\[7pt]
x&=&x_1+(z-z_1)\delta_2  \,,
\end{array}
\end{equation}
\noindent where $\delta_1=(o(|y_1|),\cdots,o(|y_1|))^T$, $\delta_2=(O(|y_1|^\alpha),\dots,O(|y_1|^\alpha))^T$ and $\alpha$ is determined by the spectrum gap between the week stable eigenvalue and the first strong stable eigenvalue (i.e. $-\lambda\pm\omega i$ and $\alpha_1$ for the system $X$ mentioned in the beginning of this paper).
\par{}
Recall the transformation (\ref{eq:pp1_3}) for the $y$-coordinate of a point on $\Pi_2$: 
\begin{equation*}
\omega\ln\dfrac{1}{|y|}=2\pi j+ \xi -\theta \, ,\quad\xi\in[0,2\pi)\,.
\end{equation*}
\noindent By the formula (\ref{eq:setting_8a}) of the map $T_2$ , the formula (\ref{eq:WsQ}) for $W^{s}_{loc}(Q_1)$, and Lemma \ref{lem:index_2} (index-2 condition), finding a periodic orbit $\{Q_1,Q_2,Q_3\}$ of period 3 and index 2 with the property $M^+\in W^s(Q_1)$ is equivalent to solve the following system of equations:
\begingroup
\allowdisplaybreaks
\begin{align*}
\quad\quad  & y_1=-\mu-Bx_3|y_3|^\rho\cos\xi_3+o(|y_3|^\rho) \,,\\[10pt]
& x_1=1+\zeta+B_1x_3|y_3|^\rho\cos(\xi_3+\theta_1-\theta)+o(|y_3|^\rho) \,,\\[10pt]
& z_1=z^- + 
                \begin{pmatrix}
                 B_2x_3|y_3|^\rho\cos{(\xi_3+\theta_2-\theta)}+o(|y_3|^\rho) \\
                 \cdots \\
                  B_{n-2}x_3|y_3|^\rho\cos(\xi_3+\theta_{n-2}-\theta)+o(|y_3|^\rho) 
                \end{pmatrix}  \,,\\[10pt]
& y_2=-\mu-Bx_1|y_1|^\rho\cos\xi_3+o(|y_1|^\rho) \,,\\[10pt]
& x_2=1+\zeta+B_1x_1|y_1|^\rho\cos(\xi_1+\theta_1-\theta)+o(|y_1|^\rho) \,,\\[10pt]
& z_2=z^- + 
                \begin{pmatrix}
                 B_2x_1|y_1|^\rho\cos{(\xi_1+\theta_2-\theta)}+o(|y_1|^\rho) \\
                 \cdots \\
                  B_{n-2}x_1|y_1|^\rho\cos(\xi_3+\theta_{n-2}-\theta)+o(|y_1|^\rho) 
                \end{pmatrix}  \,,\\[10pt]
& y_3=-\mu-Bx_2|y_2|^\rho\cos\xi_2+o(|y_2|^\rho) \,,\\[10pt]
& x_3=1+\zeta+B_1x_2|y_2|^\rho\cos(\xi_2+\theta_1-\theta)+o(|y_2|^\rho) \,,\\[10pt]
& z_3=z^- + 
                \begin{pmatrix}
                 B_2x_2|y_2|^\rho\cos{(\xi_2+\theta_2-\theta)}+o(|y_2|^\rho) \\
                 \cdots \\
                  B_{n-2}x_2|y_2|^\rho\cos(\xi_2+\theta_{n-2}-\theta)+o(|y_2|^\rho) 
                \end{pmatrix}  \,,\\[10pt]
& \mu=y_1+\delta_1(z^+-z_1)\,,\\[10pt]                
& 1=x_1+\delta_2(z^+-z_1)\,,\\[10pt]
& \cos(\xi_1-\varphi)\cos(\xi_2-\varphi)\cos(\xi_3-\varphi)=c\psi\,,
\end{align*}
\endgroup
where the first nine equations give us a periodic orbit of period 3, the next two equations imply $W^{ss}(Q)\cap M^+\neq\emptyset$, and the last one makes this orbit having index 2. After expressing $x$ and $z$ as functions of $y$, we can drop the equations for them (except the one for $x_1$ used to obtain $M^+\in W^s(Q_1)$). The reduced system assumes the form
\begingroup
\allowdisplaybreaks
\begin{align} 
\quad\quad  & y_1=-\mu-B_{\mbox{new}}|y_3|^\rho\cos\xi_3+o(|y_3|^\rho)
+ O(|y|_2^\rho |y|_3^\rho) +O(|y|_1^\rho |y|_2^\rho |y|_3^\rho) \,, \label{eq:pp2_1}\\[10pt] 
& y_2=-\mu-B_{\mbox{new}}|y_1|^\rho\cos\xi_1+o(|y_1|^\rho) 
+ O(|y|_1^\rho |y|_3^\rho) +O(|y|_1^\rho |y|_2^\rho |y|_3^\rho)\,, \label{eq:pp2_2}\\[10pt]
& y_3=-\mu-B_{\mbox{new}}|y_2|^\rho\cos\xi_2+o(|y_2|^\rho) 
+ O(|y|_1^\rho |y|_2^\rho) +O(|y|_1^\rho |y|_2^\rho |y|_3^\rho)\,,  \label{eq:pp2_3}\\[10pt]
& x_1=1+\zeta+B_1|y_3|^\rho\cos(\xi_3+\theta_1-\theta)+o(|y_3|^\rho) \,,\label{eq:pp2_4}\\[10pt] 
& \mu=y_1+\delta_1(z^+-z_1)\,,\label{eq:pp2_5}  \\[10pt]              
& 1=x_1+\delta_2(z^+-z_1)\,,\label{eq:pp2_6}\\[10pt]  
& \cos(\xi_1-\varphi)\cos(\xi_2-\varphi)\cos(\xi_3-\varphi)=c\psi\,, \label{eq:pp2_7}
\end{align}
\endgroup
\noindent where $B_{\mbox{new}}=(1+\zeta)B$ and we drop the subscript for simplicity. We now impose two relations among $y_1,y_2$ and $y_3$ which are 
\begin{equation}\label{eq:pp2_10}
|y_3|^\rho \sim |y_1| \quad \mbox{and} \quad |y_1|^\rho \sim |y_2|.
\end{equation}
\noindent It can be seen later that these relations agree with the solutions to above system of equations. Therefore, we replace the last three terms in each of equations (\ref{eq:pp2_1}) - (\ref{eq:pp2_3}) by $o(|y_3|^\rho)$, $o(|y_1|^\rho)$ and $o(|y_2|^\rho)$, respectively. 
\par{}
From now on, we will denote by dots the small terms which tend to zero as $j_1,j_2,j_3\to +\infty$. By plugging equation (\ref{eq:pp2_6}) into (\ref{eq:pp2_4}) and letting 
\begin{equation}\label{eq:zeta}
\zeta=-\delta_2(z^+-z_1)\,,
\end{equation}
\noindent we have
\begin{equation}\label{eq:pp2_7.1}
B_1|y_3|^\rho\cos(\xi_3+\theta_1-\theta)+o(|y_3|^\rho)=0\,,
\end{equation}
\noindent which, by dividing $|y_3|^\rho$ on both sides, gives
\begin{equation}\label{eq:pp2_8}
B_1\cos(\xi_3+\theta_1-\theta)+o(1)_{y_3\to 0}=0\,.
\end{equation}
\noindent This implies 
\begin{equation}\label{eq:pp2_9}
\xi_3=\dfrac{\pi}{2}+k\pi-\theta_1+\theta+...\quad \mbox{and}\quad \cos{\xi_3}=\sin(\theta_1-\theta)+...\,.
\end{equation}
\noindent Note that, to obtain (\ref{eq:pp2_8}), we only need $\zeta+\delta_2(z^+-z_1)\sim o(|y_3|^\rho)$. Recall the assumption at the beginning of the proof that $(\theta_1-\theta)\in(0,\dfrac{\pi}{2}]$. We have $k=0$ by noting $\xi_3\in[0,2\pi)$. We remark here that the relation $|y_3|^\rho \sim |y_1|$ can now be obtained by plugging (\ref{eq:pp2_5}) and (\ref{eq:pp2_9}) into (\ref{eq:pp2_1}), which is one of the relations \ref{eq:pp2_10} we assumed before.
\par{}
By using the relations given by (\ref{eq:pp2_10}) and plugging equation (\ref{eq:pp2_5}) into (\ref{eq:pp2_1}) - (\ref{eq:pp2_3}), we get
\begingroup
\allowdisplaybreaks
\begin{align} 
\quad\quad  & y_1=-y_1-B|y_3|^\rho\cos\xi_3+o(|y_3|^\rho) \,, \label{eq:pp2_11}\\[10pt] 
& y_2=-y_1-B|y_1|^\rho\cos\xi_1+o(|y_1|^\rho) \,, \label{eq:pp2_12}\\[10pt]
& y_3=-y_1-B|y_2|^\rho\cos\xi_2+o(|y_2|^\rho) \,.  \label{eq:pp2_13}
\end{align}
\endgroup
\noindent 
By plugging equation (\ref{eq:pp2_9}) into (\ref{eq:pp2_11}), we further obtain
\begin{equation}\label{eq:pp2_14}
y_1=-\frac{1}{2}B|y_3|^\rho \sin(\theta_1-\theta)+o(|y_3^\rho|)\,.
\end{equation}  
We now apply the transformation (\ref{eq:pp1_3}) to equations (\ref{eq:pp2_14}), (\ref{eq:pp2_12}) and (\ref{eq:pp2_13}), and obtain 
\begingroup
\allowdisplaybreaks
\begin{align} 
\quad
&\quad\mathrm{exp}\bigg(\dfrac{-2\pi j_1-\xi_1+\theta}{\omega}\bigg)  =
\frac{1}{2}B\mathrm{exp}\bigg(\dfrac{(-2\pi j_3-\xi_3+\theta)\rho}{\omega}\bigg)\sin{\sigma}
+o\bigg(\mathrm{exp}\bigg(\dfrac{-2\pi \rho j_3}{\omega}\bigg)\bigg)\,, \label{eq:pp2_15} \\[30pt]
&-\mathrm{exp}\bigg(\dfrac{-2\pi j_2-\xi_2+\theta}{\omega}\bigg) =
\mathrm{exp}\bigg(\dfrac{-2\pi j_1-\xi_1+\theta}{\omega}\bigg)
-B\mathrm{exp}\bigg(\dfrac{(-2\pi j_1-\xi_1+\theta)\rho}{\omega}\bigg)\cos{\xi_1} \notag
\\[20pt]
&\qquad\qquad\qquad\qquad\qquad\qquad+o\bigg(\mathrm{exp}\bigg(\dfrac{-2\pi \rho j_1}{\omega}\bigg)\bigg)\,, \label{eq:pp2_16} \\[30pt]
&-\mathrm{exp}\bigg(\dfrac{-2\pi j_3-\xi_3+\theta}{\omega}\bigg ) =
\mathrm{exp}\bigg(\dfrac{-2\pi j_1-\xi_1+\theta}{\omega}\bigg)
-B\mathrm{exp}\bigg(\dfrac{(-2\pi j_2-\xi_2+\theta)\rho}{\omega}\bigg)\cos{\xi_2} \notag
\\[20pt]
&\qquad\qquad\qquad\qquad\qquad\qquad+o\bigg(\mathrm{exp}\bigg(\dfrac{-2\pi \rho j_2}{\omega}\bigg)\bigg)\,. \label{eq:pp2_17}
\end{align}
\endgroup
\noindent Let us solve those equations for sufficiently large $j_1,j_2$ and $j_3$. 
\par{}
We first divide equation (\ref{eq:pp2_15}) by $o\big(\mathrm{exp}\big(\dfrac{-2\pi \rho j_3}{\omega}\big)\big)$ on both sides,  and then take $j_3$ large enough. After taking logarithm on both sides of the resulting equation, we obtain
\begin{equation}\label{eq:pp2_18}
\rho j_3- j_1=\frac{1}{2\pi}
\big(\omega\log{\big(\frac{1}{2}B\sin{\sigma}\big)}-\theta+\rho\theta-\rho\xi_3+\xi_1\big)+\ldots \,.
\end{equation}
\noindent In a similar way, equation (\ref{eq:pp2_16}) gives 
\begin{equation}\label{eq:pp2_19}
\cos{\xi_1}=B^{-1}\mathrm{exp}\bigg(\dfrac{2\pi(\rho j_1- j_2)}{\omega}\bigg)
\mathrm{exp}\bigg(\dfrac{\theta-\rho \theta+\rho \xi_1-\xi_2}{\omega}\bigg)
+\ldots \, .
\end{equation}
\noindent By moving the first term on the RHS of (\ref{eq:pp2_17}) to its LHS, multiplying $-1$ on both sides, and then taking logarithm, we have
\begin{equation}\label{eq:pp2_21}
\log{\bigg(\mathrm{exp}\bigg(\dfrac{-2\pi j_3-\xi_3+\theta}{\omega}\bigg ) +
\mathrm{exp}\bigg(\dfrac{-2\pi j_1-\xi_1+\theta}{\omega}\bigg)\bigg)}
=\log{\bigg(B\mathrm{exp}\bigg(\dfrac{(-2\pi j_2-\xi_2+\theta)\rho}{\omega}\bigg)\cos{\xi_2}+\cdots\bigg)}\,, \\[20pt]
\end{equation}
\noindent i.e.
\begin{equation}\label{eq:pp2_22}
\dfrac{-2\pi j_1-\xi_1+\theta}{\omega}+\log\bigg(1+\mathrm{exp}\bigg(\dfrac{-2\pi(j_3-j_1)-\xi_3+\xi_1}{\omega}\bigg)
=\log{\bigg(B\mathrm{exp}\bigg(\dfrac{(-2\pi j_2-\xi_2+\theta)\rho}{\omega}\bigg)\cos{\xi_2}+\cdots\bigg)}\,.
\end{equation}
\noindent By noting $j_3\gg j_1$ from the relation $|y_3|^\rho \sim |y_1|$ stated in (\ref{eq:pp2_10}), the last equation implies
\begin{equation}\label{eq:pp2_20}
\cos{\xi_2}=B^{-1}\mathrm{exp}\bigg(\dfrac{2\pi(\rho j_2- j_1)}{\omega}\bigg)
\mathrm{exp}\bigg(\dfrac{\theta-\rho \theta+\rho \xi_2-\xi_1}{\omega}\bigg)
+\ldots \, .
\end{equation}
\par{}
Let us now look into equation (\ref{eq:pp2_7}) of the index-2 condition. The equation (\ref{eq:pp2_9}) implies that $\cos{(\xi_3-\varphi)}$ cannot be generically arbitrarily small. We show that $\cos(\xi_2-\varphi)$ also cannot be arbitrarily small. Suppose $\cos(\xi_2-\varphi)=o(1)_{j_1,j_2,j_3\to \infty}$ and note that $y_i$ are close to 0. Since $\cos\xi_2$ is finite and $|y_3|^\rho \sim |y_1|$ (\ref{eq:pp2_10}), we have $|y_2|^\rho \sim |y_1|$ from (\ref{eq:pp2_13}). This contradicts with our assumption that $|y_1|^\rho \sim |y_2|$ (\ref{eq:pp2_10}). We now assume that $\cos{(\xi_1-\varphi)}=o(1)_{j_1,j_2,j_3\to \infty}$, which leads to
\begin{equation}\label{eq:pp2_23}
\xi_1=\frac{3\pi}{2}+\varphi+\ldots  \quad \mbox{and} \quad\cos{\xi_1}=\sin{\varphi}+\ldots \quad .
\end{equation}
We plug equation (\ref{eq:pp2_23}) into (\ref{eq:pp2_19}), and then get
\begin{equation}\label{eq:pp2_24}
\rho j_1- j_2=\frac{1}{2\pi}
\big(\omega\log{(B\cos{\xi_1})}-\theta+\rho\theta-\rho\xi_1+\xi_2\big)+\ldots \quad.
\end{equation}
\noindent Since $j_1,j_2$ are large and the RHS of (\ref{eq:pp2_24}) is uniformly bounded, we have\\[-5pt] \[\rho j_1 \sim j_2\,,\] which agrees with the assumption that $|y_1|^\rho \sim |y_2|$ (\ref{eq:pp2_10}).
We can find $\xi_2$ by plugging this into equation (\ref{eq:pp2_20}):
\begin{equation}\label{eq:pp2_25}
\cos{\xi_2}=o(1)_{j_1,j_2,j_3\to +\infty} \quad \mbox{and} \quad \xi_2=\dfrac{\pi}{2}+\ldots \quad \text{or}\quad \dfrac{3\pi}{2}+\ldots  \quad.
\end{equation}
\noindent For certainty, we let $\xi_2=\dfrac{\pi}{2}+\cdots$. We now rewrite (\ref{eq:pp2_24}) and (\ref{eq:pp2_18}) with values of $\xi_i$ as
\begin{equation}\label{eq:pp2_26}
\rho j_1- j_2=\frac{1}{2\pi}
\bigg(\omega\log{(B\sin{\varphi})}-\theta+\rho\theta-\rho\bigg(\frac{3\pi}{2}+\varphi\bigg)+\frac{\pi}{2}\bigg)+\ldots =:u+\ldots \,,
\end{equation}
\noindent and
\begin{equation}\label{eq:pp2_27}
\rho j_3- j_1=\frac{1}{2\pi}
\bigg(\omega\log{\bigg(\frac{1}{2}B\sin{(\theta_1-\theta)}\bigg)}-\theta+\rho\theta-\rho\bigg(\frac{\pi}{2}+\theta-\theta_1\bigg)+\frac{3\pi}{2}+\varphi\bigg)+\ldots =:v+\ldots \,.
\end{equation}
Equations (\ref{eq:pp2_26}) and (\ref{eq:pp2_27}) are relations among parameters. If we can find integers $j_1,j_2$ and $j_3$ such that the parameters satisfy the these two relations, then the system of equations (\ref{eq:pp2_1}) - (\ref{eq:pp2_7}) can be solved. In fact, for any given $N$, we need (\ref{eq:pp2_26}) and (\ref{eq:pp2_27}) to be satisfied with some $(j_1,j_2,j_3)$ where $j_i>N$ $(i=1,2,3)$. This is because we need $j_i$ to be sufficiently large so that the terms denoted by dots can be sufficiently small when we take the limit $j_i \to +\infty$. 
\par{}
Now recall the parameter values of $\rho,u$ and $v$ stated in Lemma \ref{lem:nontransverse2}, which are 
\begin{equation}\label{eq:rhouv}
\hat{\rho}=\frac{p}{q}\, ,\hat{u}=\frac{p_1}{q}\, ,\hat{v}=\frac{p_2}{q}\, ,
\end{equation}
\begin{spacing}{1.3}\noindent where $p,q$ are co-prime integers and $p_1,p_2$ are any integers. We now show that there exists a sequence $\{(j_1^i,j_2^i,j_3^i)\}$ of triples of integers where $j_1^i,j_2^i,j_3^i\to +\infty$ as $i\to +\infty$ such that, for each triple $(j_1^i,j_2^i,j_3^i)$, the corresponding parameter values $u_i,v_i$ obtained from the relations (\ref{eq:pp2_26}) and (\ref{eq:pp2_27}) with $\rho=\hat{\rho}$ satisfy that $u_i\to \hat{u}$ and $v_i \to \hat{v}$ as $i\to +\infty$. Finding such sequence $\{(j_1^i,j_2^i,j_3^i)\}$ is equivalent to seeking for integer solutions $(j_1^i,j_2^i,j_3^i)$ with $j_1^i,j_2^i,j_3^i\to +\infty$ as $i\to +\infty$ to the following system of equations:\end{spacing}
\begin{equation}\label{eq:pp2_28}
\left\{\begin{array}{rcl}
      \hat{\rho} j_1 -j_2&=&\hat{u} \,,\\[5pt]
      \hat{\rho} j_3 -j_1&=&\hat{v} \,.
\end{array}\right.
\end{equation} 
\noindent By plugging the relations (\ref{eq:rhouv}) into (\ref{eq:pp2_28}), we get a system of two linear Diophantine equations:
\begin{equation}\label{eq:pp2_29}
\left\{\begin{array}{rcl}
      p j_1 -q j_2&=&p_1\,, \\[5pt]
      p j_3 -q j_1&=&p_2\,.
\end{array}\right.\
\end{equation} 
\begin{spacing}{1.2}
\noindent Note that a linear Diophantine equation $ax+by=c$ has integer solutions if and only if $c$ is a multiple of $\gcd(a,b)$. For a known solution $(x,y)$, we can construct infinitely many solutions of the form $(x+kv,y-ku)$, where $u=\frac{a}{gcd(a,b)}, v=\frac{b}{gcd(a,b)}$ and $k=0,\pm1,\pm2\ldots \,.$ It is obvious that if $p, q$ are co-prime, then the two equations in (\ref{eq:pp2_29}) can be solved separately. The solutions to the first equation are of the form $(\hat{j}_1+kq, \hat{j}_2+kp)$, where $(\hat{j}_1,\hat{j}_2)$ is a solution to the first equation and $k$ is an arbitrary integer. Now let us plug $j_1=\hat{j}_1+kq$ into the second equation of (\ref{eq:pp2_29}) and sort the terms. We have\end{spacing}
\begin{equation}\label{eq:pp2_30}
pj_3-q^2k=n+q \hat{j}_1\,.
\end{equation}
\noindent Consider $j_3$ and $k$ as unknowns. Note that $p, q^2$ are co-prime since $p, q$ are co-prime. Thus, the solutions to (\ref{eq:pp2_30}) are of the form $(\hat{j}_3+iq^2, \hat{k}+ip)$, where $(\hat{j}_3,\hat{k})$ is a special solution to this equation and $i$ is an arbitrary integer. Therefore, we have infinitely many solutions $(j^i_1,j^i_2,j^i_3)=(\hat{j}_1+(\hat{k}+ip)q, \hat{j}_2+(\hat{k}+ip)p, \hat{j}_3+iq^2)$ to (\ref{eq:pp2_29}). Obviously, the integers $j^i_1,j^i_2$ and $j^i_3$ can be simultaneously arbitrarily large. Hence, we find the desired sequence $\{(j_1^i,j_2^i,j_3^i)\}$.
\par{}
For each triple $(j_1^i,j_2^i,j_3^i)$ with sufficiently large $i$, the system of equations (\ref{eq:pp2_1}) - (\ref{eq:pp2_7}) can be solved. The corresponding parameter values $\mu_i$ and $\zeta_i$ are obtained from equations (\ref{eq:pp2_5}) and (\ref{eq:zeta}), respectively. Lemma \ref{lem:nontransverse2} is proved. \qed
\subsubsection{Quasi-transverse intersection $W^s(Q)\cap W^u(P)$}\label{prf:quasitrans2}
One can check, by iterating a vertical line connecting $P$ and a point in $\Pi_0$ like what we did in the proof of Theorem \ref{thm1}, that the unstable manifold of the index-1 fixed point $P(y_p,x_p,z_p)$ is a spiral winding onto $M^+(\mu,1,z^+)$, which is given by
\begin{equation}\label{eq:thm2_1}
\begin{array}{rcl}
y&=&\mu+Ax_pt^{\rho}\cos(\omega\ln(\dfrac{1}{t})+\eta)+o(t^\rho) \,,\\[7pt]
x&=&1+A_1x_pt^{\rho}\cos(\omega\ln(\dfrac{1}{t})+\eta_1)+o(t^\rho) \,, \\[7pt]
z&=&z^+ +
                \begin{pmatrix}
                 A_2x_p t^\rho\cos{(\omega \ln {\dfrac{1}{t}}+\eta_2)}+o(t^\rho) \\
                 \cdots \\
                  A_{n-2}x_p t^\rho\cos(\omega \ln {\dfrac{1}{t}}+\eta_{n-2})+o(t^\rho) 
                \end{pmatrix} \,,
\end{array}
\end{equation}
\noindent where $t\in(0,y_p)$. We take a point $Q$ given by Lemma \ref{lem:nontransverse2} at parameter values $(\mu_j,\zeta_j,\rho^*,$ $u_{j}, v_{j})$ with $j$ sufficiently large such that $P$ remains a saddle fixed point. It follows that the non-empty intersection $W^u(P)\cap W^s(Q)$ can be created by an arbitrary perturbation in $\zeta$ in system $X_{\mu_j,\zeta_j,\rho^*,u_{j}, v_{j}}$. Indeed, by changing $\zeta$, one can change the distance corresponding to $x$-coordinate between $M^+$ and $W^s_{loc}(Q_1)$ (see Figure \ref{fig:3d_2} and equation (\ref{eq:zeta})). Therefore, for each $j$, one can find a sequence $\zeta_j^i$ such that $W^u(P)\cap W^s(Q)$ is non-empty in the system $X_{\mu_j,\zeta_j^i,\rho^*,u_{j}, v_{j}}$, where $\zeta_j^i\to\zeta_j$ as $i\to +\infty$. Consequently, one can construct a new sequence $\{(\mu_{j},\zeta_j,u_{j}, v_{j})\}$ such that system $X_{\mu_j,\zeta_j,\rho^*,u_{j}, v_{j}}$ has the intersection $W^u(P)\cap W^s(Q)$. This intersection is quasi-transverse by the same argument used in the beginning of Section \ref{quasitrans1}.
\subsubsection{Transverse intersection $W^u(Q)\cap W^s(P)$}
We now prove that the point $P\in\{P^+_k\}$ given by Lemma \ref{lem:WsP_2} and the point $Q$ given by Lemma \ref{lem:nontransverse2} also satisfy $W^u(Q) \cap W^s(P)\neq\emptyset$. Specifically, we first fix the point $Q$ and then find the point $P$ by Lemma \ref{lem:WsP_2} with the $\mu$ value corresponding to $Q$.
\par{}
We remark here that we only show the existence of the transverse intersection $W^u(Q) \cap W^s(P)$ for the case where $\mu_j<0$ (the $\mu$ value corresponding to $Q$) and $(\theta_1-\theta)\in(0,\dfrac{\pi}{2}]$ i.e. the fourth configuration shown in figure \ref{fig:config}. One can easily check that the transverse intersection $W^u(Q) \cap W^s(P)$ exists in other cases as well. 
\par{}
Let $Q$ has orbit $\{Q=Q_1,Q_2,Q_3\}$. We take an arbitrarily small neighborhood $U$ of the point $Q_1$. We claim that there exists some $i$ such that $T_2^{(i)}(U)$ intersects $\Pi_0$ transversely. Indeed, this can be achieved by applying the same argument used in the proof of Lemma \ref{lem:WuQ}. We remark here that although Lemma \ref{lem:WuQ} is for the case where $\mu=0$, it also holds for small $\mu\neq0$. Indeed, the key step in the proof of Lemma \ref{lem:WuQ} is to show that the orbit of the point $\tilde{Q_1}$ (projection of the index-2 point along the a leaf of the foliation $\mathcal{F}_1$) is inside of a region in $\tilde{\Pi}$ bounded by $\{x=1+\delta\}$, $\{x=1-\delta\}$ and two pre-images of $\{y=0\}$ under $\tilde{T}$. In the case where $\mu\neq 0$, the difference is that there are only finite pre-images of $\{y=0\}$ (the smaller the value of $\mu$ is the more pre-images we get). However, there are always some pre-images have a finite distance to $\{y=0\}$ as long as $\mu$ is not too large. This implies that we can still find the desired region that contains the orbit of $\tilde{Q}_1$ by choosing $\{Q_1,Q_2,Q_3\}$ sufficiently close to $\Pi_0$ (i.e. taking $j$ sufficiently large in Lemma \ref{lem:nontransverse2}).
\par{}
Let $T_2^{(i)}(U)$ be the first iterate which intersects $\Pi_0$ transversely. We have three cases depending on which point $Q_1$, $Q_2$ or $Q_3$ is contained in $T_2^{(i)}(U)$. 
\par{}
If we have $Q_3\in T_2^{(i)}(U)$, i.e. $T_2^{(i)}(Q_1)=Q_3$, then there exists a connected component $l\subset(T_2^{(i)}(U)\cap\Pi_2)$ joining $Q_3$ and a point $M\in\Pi_0$. It follows that $T_2(l)\subset T_2^{(i+1)}(U)$ is a connected component joining  $T_2(Q_3)=Q_1\in\Pi_2$ and $T_2(M)=M^-\in\Pi_1$ since we assumed $\mu_j<0$. By Lemma \ref{lem:WsP_2}, the local stable manifold $W^s_{loc}(P)$ is a surface between $\{y=0\}$ and $\{y=|\mu_j|\}$. It follows that $T_2^{(i+1)}(U)\cap W^s_{loc}(P)\neq\emptyset$, which implies that $W^u_{loc}(Q_1)\cap W^s_{loc}(P)\neq\emptyset$. 
\begin{figure}[!h]
\centering
\includegraphics[width=0.65\textwidth]{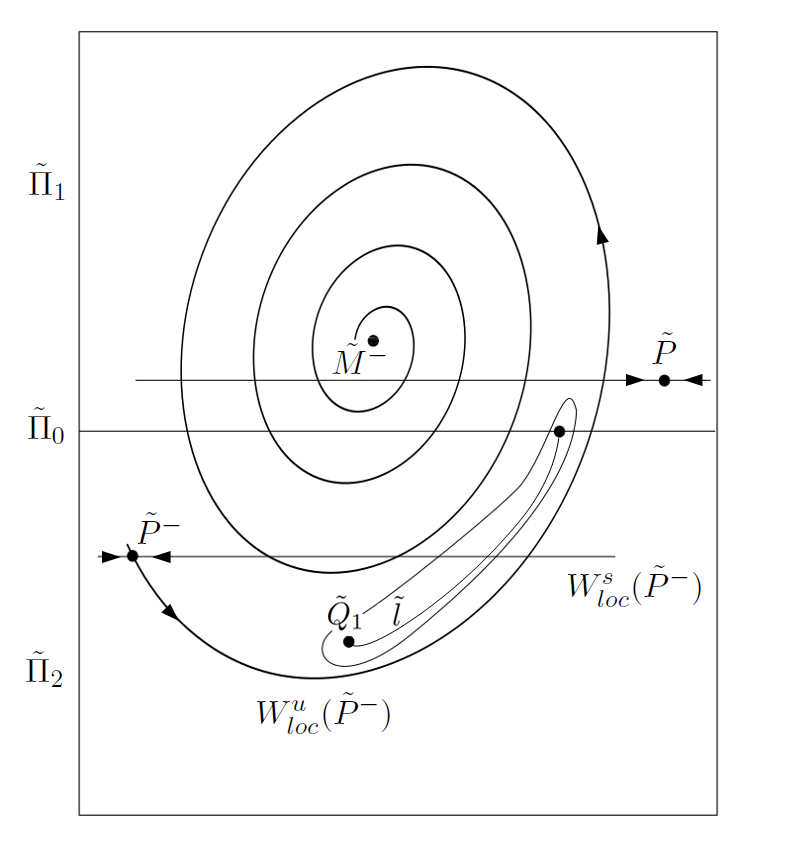}
\caption{The configuration of corresponding objects in the quotient cross-section $\tilde{\Pi}$. If $\tilde{T}_2^{(i)}(\tilde{Q}_1)=\tilde{Q}_1$ for some $i$, then we have $W^s_{loc}({\tilde{P}^-})\cap l\neq\emptyset$ and $W^u_{loc}({\tilde{P}^-})\cap W^s_{loc}({\tilde{P}})\neq\emptyset$}
\label{fig:transintersection4}
\end{figure}
\par{}
Let now $T_2^{(i)}(Q_1)=Q_1$ and $l$ be the connected component joining $Q_1$ and a point in $\Pi_0$. Note that we have $y_1\sim\mu_j$ since $M^+\in W^s(Q)$ (see (\ref{eq:pp2_5})). By applying Lemma \ref{lem:WsP_2} to the set $\{P^-_k\}$, one can find a point $P^-\in\{P^-_k\}$ such that $P^-$ remains an index-1 fixed point at $\mu=\mu_j$, and its local stable manifold $W^s_{loc}(P^-)$ intersects $l$. This gives $W^s(P^-)\cap W^u(Q_1)\neq \emptyset$. Since $W^u(P^-)$ is a spiral winding onto $M^-\in\Pi_1$, it must intersect the surface $\{y=-\mu_j\}\subset\Pi_1$. Also by Lemma \ref{lem:WsP_2}, one can find a point $P\in\{P^+_k\}$ that remains fixed at $\mu=\mu_j$ with a local stable manifold between $\{y=0\}$ and $\{y=|\mu_j|\}$.  Hence, we have $W^u(P^-)\cap W^s(P)\neq\emptyset$ (see figure \ref{fig:transintersection4}), which further implies $W^u(Q_1)\cap W^s(P)\neq\emptyset$. The same result holds if $T_2^{(i)}(Q_1)=Q_2$. Indeed, the relation $|y_1|^\rho\sim |y_2|$ implies that $|y_1|\ll |y_2|$, and, therefore, $W^s_{loc}(P^-)$ must intersect the connected component joining $Q_2$ and a point in $\Pi_0$.
\par{}
Thus, for each quadruple $(\mu_{j},\zeta_j,u_{j}, v_{j})$, the system $X_{\mu_j,\zeta_j,\hat{\rho},u_{j}, v_{j}}$ has a heterodimensional cycle related to two periodic orbits of index 2 and index 3 which correspond to the periodic points $P$ and $Q$ of the map $T$, respectively. Theorem \ref{thm2} is proved.
\par{}
We remark here that any point $P'$ from $\{P^+_k\}$ which, under the perturbation, remains an index-1 fixed point and is homoclinically related to the point $P$ given by Lemma \ref{lem:WsP_2} gives an transverse intersection $W^s(P')\cap W^u(Q)$. The quasi-transverse intersection $W^s(Q)\cap W^u(P')$ can be obtained by the same argument in Section \ref{prf:quasitrans2}. Therefore, such point $P'$ can also be used to create a heterodimensional cycle with the point $Q$.
\section{Discussion}
In this paper, we have showed two mechanisms for the emergences of heterodimensional cycles near a pair of Shilnikov loops. The further research is to check whether those heterodimensional cycles are robust or there are robust cycles nearby. Ideally, One might be able to find a blender or its analogue near a pair of Shilnikov loops.
\par{}
This also links to the application question: how to detect heterodimensional cycles? As we know that Shilnikov loop along is a simple criterion for chaos. Now a pair of such loops, provided the volume hyperbolicity and the non-coincidence condition (which are also reasonably easy to verify), gives a simple and practical criterion for the heterodimensional chaos (i.e. chaotic dynamics where saddles with different dimensions of unstable manifolds coexist and are connected).
\section*{Acknowledgement}
The author is grateful to his scientific adviser Dmitry Turaev for setting this problem and useful discussions.


\begin{thebibliography}{1}
\bibitem[An67]{an67} Anosov, D. V.. Geodesic flows on closed Riemannian manifolds of negative curvature.
Proc. Steklov Inst. Math., 90 (1967), 1–235.
%
\bibitem[ABS83]{abs} Afraimovich, V. S., Bykov, V. V. and Shilnikov, L. P.. On the structurally unstable attracting limit sets of Lorenz attractor type. Tran. Moscow Math. Soc., 2 (1983), 153-215. 

\bibitem[BD96]{bd96} Bonatti, C. and Díaz, L. J.. Persistent transitive diffeomorphisms. Annals of Mathematics, 143(2) (1996), 357-396.
\bibitem[BDV00]{bdv} Bonatti, C., Díaz, L. J. and Viana, M.. Dynamics Beyond Uniform Hyperbolicity. Springer, Berlin, Heidelberg, New York, 2000.
\bibitem[BD08]{bd08} Bonatti, C. and Díaz, L. J.. Robust heterodimensional cycles and C1-generic dynamics. J. Inst. Math. Jussieu 7, no. 3 (2008), 469-525.
\bibitem[BC15]{bc15} Bonatti, C. and Crovisier, S.. Center manifolds for partially hyperbolic sets without strong unstable connections. Journal of the Institute of Mathematics of Jussieu, available on CJO2015. doi:10.1017/S1474748015000055. 
\bibitem[CA10]{ca10} Chawanya, T. and Ashwin, P.. A minimal system with a depth-two heteroclinic network. Dyn. Syst., 25 (2010), pp. 397–412. 
%
\bibitem[DGSY27]{dgsy27} Dawson, S., Grebogi, C., Sauer, T. and Yorke, A.. Obstructions to shadowing when a Lyapunov exponent fluctuates about zero. Phys. Rev. Lett. 73, 1927 – Published 3 October 1994.
%
\bibitem[Dí92]{d92} Díaz, L. J. and Rocha, J.. Non-connected heterodimensional cycles: bifurcation and stability. Nonlinearity, 5 (1992), 1315-1341.
\bibitem[Dí95a]{d95} Díaz, L. J.. Robust nonhyperbolic dynamics and heterodimensional cycles. Ergodic Theory and Dynamical Systems, 15 (1995), 291-315.
\bibitem[Dí95b]{d95_2} Díaz, L. J.. Persistence of cycles and nonhyperbolic dynamics at the unfolding of heteroclinic bifurcations. Nonlinearity, 8 (1995), 693-715.
\bibitem[EKTS89]{ekts89} ELeonsky, V. M ., Kulagin, N. E., Turaev, D. V. and Shilnikov, L. P.. On the classification of selflocalized states of electromagnetic field within nonlinear medium. Proceedings of the iv international workshop on nonlinear and turbulent processes in physics, volume 2 (1989), 235-238. 
%
\bibitem[Ga83]{g83} Gaspard, P.. Generation of a countable set of homoclinic flows through bifurcation, Physics Letters A, Volume 97, Issues 1-2, 8 August 1983, 1-4.
%
\bibitem[GST09]{} Gonchenko, S. V., Shilnikov, L. P. and Turaev, D. V.. On global bifurcations in three-dimensional diffeomorphisms leading to wild Lorenz-like attractors. Regular and Chaotic Dynamics, 14:1 (2009), 137.
%
\bibitem[Ho96]{ho96} Homburg, A. J.. Global aspects of homoclinic bifurcations of vector fields, Mem. Amer. Math. Soc. 121 (1996), no. 578, viii+128.
%
\bibitem[LT15]{lt15} Li, D. and Turaev, D. V.. Existence of heterodimensional cycles near Shilnikov loops. Preprint. arXiv:1512.01280 [math.DS].
%
\bibitem[OS87]{os87} Ovsyannikov, I. M. and Shilnikov, L. P.. On systems with a saddle-focus homoclinic curve. Sbornik: Mathematics 58 (2) (1987) 557-574.
%
\bibitem[OS92]{os92} Ovsyannikov, I. M. and Shilnikov, L. P.. Systems with a homoclinic curve of multidimensional saddle-focus type, and spiral chaos. Math. USSR Sbornik, 73 (1992),  415-443.
%
\bibitem[Sh65]{sh65} Shilnikov, L. P.. A case of the existence of a countable number of periodic motions (Point mapping proof of existence theorem showing neighborhood of trajectory which departs from and returns to saddle-point focus contains denumerable set of periodic motions). SOVIET MATHEMATICS 6 (1965), 163-166.
%
\bibitem[Sh67]{sh67} Shilnikov, L. P.. On a Poincaré-Birkhoff problem. Sbornik: Mathematics 3 (3), 353-371.
%
\bibitem[Sh70]{sh70} Shilnikov, L. P.. A contribution to the problem of the structure of an extended neighborhood of a rough equilibrium state of saddle-focus type. Sbornik: Mathematics 10 (1) (1970), 91-102.
%
\bibitem[SSTC01]{sstc1} Shilnikov, L. P., Shilnikov, A. L., Turaev, D. V. and Chua, L. O.. Methods Of Qualitative Theory In Nonlinear Dynamics (Part I). World Sci.-Singapore, New Jersey, London, Hong Kong, 2001.
%
\bibitem[SSTC02]{sstc2} Shilnikov, L. P., Shilnikov, A. L., Turaev, D. V. and Chua, L. O.. Methods Of Qualitative Theory In Nonlinear Dynamics (Part II). World Sci.-Singapore, New Jersey, London, Hong Kong, 2001.
%
\bibitem[ST97]{st97} Shilnikov, L. P. and Turaev, D. V.. Superhomoclinic orbits and multi-pulse homoclinic loops in Hamiltonian systems with discrete symmetries, Regular and Chaotic Dynamics 2, No. 3/4 (1997), 126–138.
%
\bibitem[ST99]{st99} Shashkov, M. V. and Turaev, D. V.. An Existence theorem of smooth nonlocal center
manifolds for systems close to a system with a homoclinic loop. J. Nonlinear Sci. Vol. 9 (1999), 525-573.
%
\bibitem[Tu96]{tu96} Turaev, D. V.. On dimension of non-local bifurcational problems, International Journal of Bifurcation and Chaos, 6(5) (1996), 919-948.
%
\bibitem[Tu01]{tu01} Turaev, D. V.. Multi-pulse homoclinic loops in systems with a smooth first integral Ergodic Theory, Analysis and Efficient Simulation of Dynamical Systems ed B Fiedler (Berlin: Springer) pp 691–716.
%
\bibitem[TS98]{ts98} Turaev, D. V. and Shilnikov, L. P.. An example of a wild strange attractor. Sbornik. Math. 189(2) (1998), 291-314.
\end{thebibliography}
\end{document}